\documentclass[11pt,reqno]{amsart}

\usepackage{amsmath,amssymb, graphicx,latexsym,amsthm,mathrsfs,amsthm}
\usepackage{graphicx,color,cite}
\usepackage{caption, subcaption}
\newtheorem{theorem}{Theorem}
\newtheorem{lemma}{Lemma}
\newtheorem{proposition}{Proposition}

\newtheorem{corollary}{Corollary}

\theoremstyle{remark}

\begin{document}

\title[Inverse scattering problems]{Increasing stability estimates for the inverse potential scattering problems}

\author[J. Zhai]{Jian Zhai}
\address{School of Mathematical Science, Fudan University, Shanghai 200433 China}
\email{jianzhai@fudan.edu.cn}
  
\author[Y. Zhao]{Yue Zhao}
\address{School of Mathematics and Statistics, Central China Normal University, Wuhan 430079, China}
\email{zhaoyueccnu@163.com}

%\thanks{}



\begin{abstract}

This paper is mainly concerned with the inverse scattering problem of determining the unknown potential for the classical Schr\"odinger equation in two and three dimensions.
We establish the increasing stability of the inverse scattering problem from either multi-frequency near-field data or multi-frequency far-field pattern. 
The stability estimate consists of the Lipschitz type data discrepancy and the logarithmic high frequency tail of the potential function, where the latter
decreases as the upper bound of the frequency increases. A novel method is proposed for the proof, which is based on choosing appropriate incident plane waves 
and an application of the quantitative analytic continuation. A key ingredient in the analysis is employing scattering theory to obtain an analytic region and resolvent
estimates in this region for the resolvent in two and three dimensions. 
We further apply this method to study the inverse scattering problem of determining both the magnetic potential and electric potential for the three-dimensional magnetic Schr\"odinger equation.
\end{abstract}

\maketitle

\section{Introduction}

We are mainly concerned with the increasing stability of determining a compactly supported unknown potential in the Schr\"odinger equation in $\mathbb R^n, n = 2, 3$. Let $e^{\mathrm{i}\kappa \boldsymbol x\cdot \boldsymbol d}$ be an incident plane wave with direction $\boldsymbol d\in\mathbb S^{n-1}$.
We consider the scattering problem modeled by the following classical Schr\"odinger equation in $\mathbb R^n$
\begin{align}\label{main_eq}
-\Delta u - \kappa^2 u+ Vu = 0, 
\end{align}
where $u(\boldsymbol x,\kappa, \boldsymbol d) = e^{\mathrm{i}\kappa \boldsymbol x\cdot \boldsymbol d} + u^s(\boldsymbol x,\kappa,\boldsymbol d)$ is the total field, $u^s$ is the scattered field,
$\kappa>0$ is the wavenumber and $V(\boldsymbol x)\in L_{\rm comp}^\infty(\mathbb R^n)$ is the potential function. We assume that the
support of $V(\boldsymbol x)$ is contained in $B_R := \{\boldsymbol x\in\mathbb R^{n} ~:~ |\boldsymbol x|\leq R\}$ where $R>0$ is a constant. 
Let $\partial B_R$ be the boundary of $B_R$. The Sommerfeld radiation condition is imposed on the scattered field to ensure the well-posedness of the
problem:
\begin{align}\label{src}
\lim_{r\to\infty}r^{\frac{n - 1}{2}} (\partial_r u^s-\mathrm{i}\kappa u^s)=0,\quad r=|\boldsymbol x|
\end{align}
uniformly in all directions $\hat{\boldsymbol x} = \boldsymbol x/|\boldsymbol x|$.
We are interested in the inverse scattering problem of determining $V$ from the near-field Dirichlet boundary measurements 
$u^s(\boldsymbol x, \kappa, \boldsymbol d)\vert_{\partial B_R}$, $\boldsymbol x\in\partial B_R$, $\boldsymbol d\in \mathbb S^{n - 1}$ or the far-field pattern corresponding to the wavenumber $\kappa$ given in a
finite interval. 

There is extensive literature on the inverse boundary value problems of determining the potential from the \textit{Dirichlet-to-Neumann} (DtN) map for the Schr\"odinger equation at a fixed frequency. For instance, the uniqueness result was proved in \cite{SU} and a logarithmic stability estimate was obtained in \cite{Ale}. Increasing stability estimates 
in three dimensions were proved in \cite{Isakov, Isakov1, Isakov2, NUW}. Also a linearized problem was studied in \cite{ILX} with numerical evidence for increasing stability.
For such inverse boundary value problems, the main tool is the construction of \textit{Complex Geometric Optics} (CGO) solutions. We point it out that so far the method of constructing CGO solutions cannot be used to obtain increasing stability estimates in two dimensions \cite{NUW}.  On the other hand, for the inverse scattering problems, since the knowledge of the DtN map is equivalent to knowing the far-field pattern at a fixed frequency \cite{Salo}, many uniqueness results follow from those of the corresponding inverse boundary value problems. 
However, compared with the inverse boundary value problems, stability estimates for the inverse scattering problems are much less studied. 
We mention that stability estimates of recovering the refractive index at a fixed frequency were obtained in \cite{HH, IN} based
on CGO solutions, which are of logarithmic type as well. This type of stability estimate reveals the ill-posed nature of the inverse problem, which
means that a small variation of the data might lead to a huge error in the reconstruction. Thus, computationally it is crucial to investigate the
stability of the inverse problem. In this work we hope to improve the stability using multi-frequency measurements in both two and three dimensions.

Recently, it has been observed experimentally (see e.g. \cite{BLLT}) that by multi-frequency measurements one may overcome the ill-posed nature of the inverse scattering problems and obtain increasing stability. Such observations have been corroborated mathematically for the inverse source problems. Specifically, the increasing stability for the inverse source problems with multi-frequency measurements has been obtained for various wave equation (see e.g. \cite{blz, CIL, IL, LY, LZZ}).
It is worth mentioning that an increasing stability estimate was proved in \cite{bt} for the multi-frequency inverse medium scattering problem in one dimension. However, the increasing stability issue for the inverse scattering problem of determining the potential function in higher dimensions remains widely open. Compared with the inverse source problem, the inverse scattering problem of determining the potential is more involved due to its nonlinearity.
In this paper, we show an increasing stability estimate using either the multi-frequency near-field data or the far-field pattern. 
The stability estimate consists of the Lipschitz type data discrepancy and the high frequency tail of the potential function, where the latter
decreases as the upper bound of the frequency increases. The proof is based on choosing appropriate incident plane waves 
and an application of the quantitative analytic continuation. We note here that, with multifrequency data available, we do not resort to the use of CGO solutions.
A key ingredient in the analysis is employing scattering theory to obtain an analytic region and resolvent
estimates in this region for the resolvent in two and three dimensions.

The method proposed in this paper can also be used to study the magnetic Schr\"odinger equation. We are interested in the determination of both the magnetic potential and electric potential. For the inverse problems, it is known that there is a natural obstruction to the uniqueness of the magnetic potential. As a matter of fact, changing the magnetic potential to its gauge equivalence would not change the measurements. For the inverse boundary value problems at a fixed frequency,
 we refer the readers to \cite{NSU, Salo06, KG} and references therein for uniqueness results, and to \cite{CP, Tzou} and references therein for stability estimates.
The main tool of these works is also the construction of CGO solutions.
We also refer the reader to \cite{ER, Salo} for uniqueness results on the inverse scattering problem of the magnetic Schr\"odinger equation. However, we are not aware of any increasing stability results so far for the recovery of the magnetic potential.

The paper is organized as follows. In Section \ref{acoustic}, we
study the stability of the inverse potential scattering problem for the classical Schr\"odinger equation in two and three dimensions.
Section \ref{three} is devoted to the case of three dimensions.
The two-dimensional case is discussed in Section \ref{two}.
The stability for the magnetic Schr\"odinger equation is derived in Section \ref{magnetic}.

\section{The classical Schr\"odinger equation}\label{acoustic}

\subsection{Three-dimensional case}\label{three}
In this section we consider the Schr\"odinger equation  \eqref{main_eq} in $\mathbb R^3$. Denote the outgoing resolvent of the elliptic operator $-\Delta - \lambda^2 +V$ by 
\begin{equation}\label{resolvent}
R_V(\lambda) = (-\Delta - \lambda^2 +V)^{-1}.
\end{equation}
First we notice that $R_V(\lambda): L^2(\mathbb R^3)\to L^2(\mathbb R^3)$ in \eqref{resolvent} is well defined for $\Im\lambda\gg 1$. It is known that $R_V$ has a meromorphic continuation to the whole complex plane as the following family of operators (cf.\cite[Theorem 3.8]{Dyatlov})
\[
R_V(\lambda): L_{\text{comp}}^2(\mathbb{R}^3)\rightarrow L^2_{\text{loc}}(\mathbb{R}^3),\quad \lambda\in\mathbb{C}.
\]
Here
\[
L_{\text{comp}}^2(\mathbb{R}^3):\{u\in L^2(\mathbb{R}^3),\exists R>0, |\boldsymbol x|>R\Rightarrow |u(\boldsymbol x)|=0\},
\]
\[
L^2_{\text{loc}}(\mathbb{R}^3):\{u\in\mathcal{D}'(\mathbb{R}^3),\forall \chi\in C_c^\infty(\mathbb{R}^3),\chi u\in L^2(\mathbb{R}^3)\}.
\]
By \cite[Theorem 3.33]{Dyatlov} $R_V(\lambda)$ has no poles for $\lambda\in\mathbb{R}\setminus\{0\}$. In particular $R_V(\kappa)$ is well defined for $\kappa>0$.
Then by equation \eqref{main} and \cite[Theorem 3.37]{Dyatlov} we have $u^s(\boldsymbol x,\kappa, \boldsymbol d)= -R_V(\kappa)(V e^{\mathrm{i}\kappa\langle \bullet, \boldsymbol d\rangle})$. Subject to the Sommerfeld radiation condition, $u^s$ admits the asymptotics
\begin{equation}\label{far_asymp}
u^s(\boldsymbol x,\kappa, \boldsymbol d)=\frac{e^{\mathrm{i}\kappa|\boldsymbol x|}}{|\boldsymbol x|}\mathcal{A}_{\infty}\left(\kappa,\frac{\boldsymbol x}{|\boldsymbol x|},\boldsymbol d\right)+\mathcal{O}(|\boldsymbol x|^{-2}).
\end{equation}
Here $\mathcal{A}_\infty(\kappa,\boldsymbol \theta, \boldsymbol d)$ is called the far-field pattern.

The inverse scattering problem considered in this section is to determine $V$ from the near-field Dirichlet boundary measurements 
$u^s(\boldsymbol x, \kappa, \boldsymbol d)\vert_{\partial B_R}$, $\boldsymbol x\in\partial B_R$, $\boldsymbol d\in \mathbb S^2$ or the far-field pattern $\mathcal{A}_\infty(\kappa,\boldsymbol \theta, \boldsymbol d)$, $(\boldsymbol \theta, \boldsymbol d)\in \mathbb S^2\times \mathbb S^2$ corresponding to the wavenumber $\kappa$ given in a
finite interval.

\subsubsection{Near-field data}\label{near}

In this section we study the increasing stability estimate using multi-frequency near-field data. %We first state the resolvent estimate given in \cite{Dyatlov} for $R_V(\kappa)$. 

Let $s>0$ be an arbitrary positive constant. Define a real-valued function space
\[
\mathcal C_Q = \{V \in H^{s}(\mathbb R^3)\cap L^\infty(\mathbb R^3):  \|V\|_{H^{s}(\mathbb R^3)}\leq Q, \|V\|_{L^{\infty}(\mathbb R^3)}\leq Q,~ \text{supp}
V\subset B_R, ~ V: \mathbb R^3 \rightarrow \mathbb R \}.
\]
Assume that $V_1$ and $V_2$ are two potential functions. Let $\kappa\in(0,+\infty)$ and $\boldsymbol d_1\in\mathbb S^2, \boldsymbol d_2 \in \mathbb S^2$.
Denote the scattered field corresponding to the incident field $e^{\mathrm{i}\kappa \boldsymbol x\cdot \boldsymbol d_j}$
and $V_i$ by 
\begin{align}\label{data}
u^s_i(\boldsymbol x, \kappa, \boldsymbol d_j) = -R_{V_i}(\kappa)(V_i e^{\mathrm{i}\kappa \boldsymbol x\cdot \boldsymbol d_j}), \quad 1\leq i, j \leq 2,
\end{align}

We will use the following properties of the outgoing resolvent $R_V$.
\begin{lemma}\textnormal{(\cite[Theorem 3.10]{Dyatlov})}\label{bound}
Fix a cutoff function $\rho\in C_0^\infty(\mathbb R^n)$ where $n\geq 3$ is odd. Suppose $V\in L^\infty_{\mathrm{comp}}(\mathbb{R}^n)$.
There exist constants $A, C, C_0$ and $T$ such that, for $j=0,1,2$,
\begin{align}
\|\rho R_V(\lambda)\rho\|_{L^2(\mathbb R^n) \to H^j(\mathbb R^n)}\leq
C|\lambda|^{j - 1} e^{T(\mathrm{Im}\lambda)_-},
\end{align}
for 
\[
\Im\lambda\geq -A-\delta\log(1+|\lambda|),\quad|\lambda|>C_0,
\]
where $\delta<1/\mathrm{diam}(\mathrm{supp} V)$.
\end{lemma}
%where $K>K_0>C_0$ and $C_0$ is as specified in Lemma \ref{ac}.
Denote
\begin{align*}
u_j=e^{\mathrm{i}\kappa \boldsymbol x\cdot \boldsymbol d_j}-R_{V_j}(\kappa)(V_j e^{\mathrm{i}\kappa \boldsymbol x\cdot \boldsymbol d_j}),\quad j=1,2,
\end{align*}
which satisfy
\begin{align}\label{eqn1}
-\Delta u_j - \kappa^2 u_j + V_j u_j = 0, \quad j = 1, 2.
\end{align}
Denote $u=u_1-u_2$, $V=V_2-V_1$. Subtracting the equation \eqref{eqn1} when $j=1$ by itself when $j = 2$ gives
\begin{align}\label{eqn2}
-\Delta u - \kappa^2 u + V_1 u = V u_2.
\end{align}
Multiplying both sides of \eqref{eqn2} by $u_1$ and integrating by parts over $B_R$ gives
\begin{equation}\label{integral_id}
\int_{B_R}Vu_1u_2\mathrm{d}x=\int_{\partial B_R} -u_1\partial_\nu u+u\partial_\nu u_1\mathrm{d}s=\int_{\partial B_R}u_1\partial_\nu u_2 -u_2\partial_\nu u_1\mathrm{d}s.
\end{equation}

Also let
\[
\begin{split}
w_1&=e^{\mathrm{i}\kappa \boldsymbol x\cdot \boldsymbol d_2}-R_{V_1}(\kappa)(V_1e^{\mathrm{i}\kappa \boldsymbol x\cdot \boldsymbol d_2}),\\
w_2&=e^{\mathrm{i}\kappa \boldsymbol x\cdot \boldsymbol d_1}-R_{V_2}(\kappa)(V_2e^{\mathrm{i}\kappa \boldsymbol x\cdot \boldsymbol d_1}).
\end{split}
\]
Denote $w=w_1-w_2$. Similar to \eqref{integral_id}, we have
\begin{equation}\label{integral_id1}
\int_{B_R}Vw_1w_2\mathrm{d}x=\int_{\partial B_R} -w_1\partial_\nu w+w\partial_\nu w_1\mathrm{d}s=\int_{\partial B_R} w_1\partial_\nu w_2-w_2\partial_\nu w_1\mathrm{d}s.
\end{equation}
We can write
\begin{align}\label{identity1}
&\int_{B_R}Vu_1u_2\mathrm{d}x+\int_{B_R}Vw_1w_2\mathrm{d}x\notag\\
%=&-\int_{\partial B_R} u_1\partial_\nu u-u\partial_\nu u_1\mathrm{d}s-\int_{\partial B_R} w_1\partial_\nu w-w\partial_\nu w_1\mathrm{d}s\\
=&\int_{\partial B_R} u_1\partial_\nu u_2-u_2\partial_\nu u_1\mathrm{d}s+\int_{\partial B_R} w_1\partial_\nu w_2-w_2\partial_\nu w_1\mathrm{d}s\notag\\
=&\int_{\partial B_R} (u_1+w_1)\partial_\nu (u_2+w_2)-(u_2+w_2)\partial_\nu (u_1+w_1)\mathrm{d}s\notag\\
&-\int_{\partial B_R} u_1\partial_\nu w_2-w_2\partial_\nu u_1\mathrm{d}s\notag\\
&+\int_{\partial B_R} u_2\partial_\nu w_1-w_1\partial_\nu u_2\mathrm{d}s.
\end{align}
Notice
\begin{align}\label{iden1}
&(u_1+w_1)\partial_\nu (u_2+w_2)-(u_2+w_2)\partial_\nu (u_1+w_1)\notag\\
=&-(u_1+w_1)\partial_\nu(u_1-w_2-u_2+w_1)+(u_1-w_2-u_2+w_1)\partial_\nu(u_1+w_1)\notag\\
=&-(u_1+w_1)\partial_\nu\Big(-R_{V_1}(\kappa)(V_1e^{\mathrm{i}\kappa \boldsymbol x\cdot \boldsymbol d_1})+R_{V_2}(\kappa)(V_2e^{\mathrm{i}\kappa \boldsymbol x\cdot \boldsymbol d_1})\notag\\
&\quad\quad\quad\quad\quad\quad\quad+R_{V_1}(\kappa)(V_1e^{\mathrm{i}\kappa \boldsymbol x\cdot \boldsymbol d_2})-R_{V_2}(\kappa)(V_2e^{\mathrm{i}\kappa \boldsymbol x\cdot \boldsymbol d_2})\Big)\notag\\
&+\Big((-R_{V_1}(\kappa)(V_1e^{\mathrm{i}\kappa \boldsymbol x\cdot \boldsymbol d_1})+R_{V_2}(\kappa)(V_2e^{\mathrm{i}\kappa \boldsymbol x\cdot \boldsymbol d_1}) \notag\\
&\quad\quad+ R_{V_1}(\kappa)(V_1e^{\mathrm{i}\kappa \boldsymbol x\cdot \boldsymbol d_2})-R_{V_2}(\kappa)(V_2e^{\mathrm{i}\kappa \boldsymbol x\cdot \boldsymbol d_2})\Big)\partial_\nu(u_1+w_1),
\end{align}
and
\begin{align}\label{iden2}
&u_1\partial_\nu w_2-w_2\partial_\nu u_1\notag\\
=&-u_1\partial_\nu(u_1-w_2)+(u_1-w_2)\partial_\nu u_1\notag\\
=&-u_1\partial_\nu \Big(-R_{V_1}(\kappa)(V_1e^{\mathrm{i}\kappa \boldsymbol x\cdot \boldsymbol d_1})+R_{V_2}(\kappa)(V_2e^{\mathrm{i}\kappa \boldsymbol x\cdot \boldsymbol d_1})\Big)\notag\\
&+\Big(-R_{V_1}(\kappa)(V_1e^{\mathrm{i}\kappa \boldsymbol x\cdot \boldsymbol d_1})+R_{V_2}(\kappa)(V_2e^{\mathrm{i}\kappa \boldsymbol x\cdot \boldsymbol d_1})\Big)\partial_\nu u_1,
\end{align}
together with
\begin{align}\label{iden3}
&u_2\partial_\nu w_1-w_1\partial_\nu u_2\notag\\
=&-u_2\partial_\nu(u_2-w_1)+(u_2-w_1)\partial_\nu u_2\notag\\
=&-u_2\partial_\nu \Big(-R_{V_1}(\kappa)(V_1e^{\mathrm{i}\kappa \boldsymbol x\cdot \boldsymbol d_2})+R_{V_2}(\kappa)(V_2e^{\mathrm{i}\kappa \boldsymbol x\cdot \boldsymbol d_2})\Big)\notag\\
&+\Big(-R_{V_1}(\kappa)(V_1e^{\mathrm{i}\kappa \boldsymbol x\cdot \boldsymbol d_1})+R_{V_2}(\kappa)(V_2e^{\mathrm{i}\kappa \boldsymbol x\cdot \boldsymbol d_1})\Big)\partial_\nu u_2.
\end{align}
By Lemma \ref{bound}, we have
\[
\|u^s_i(\boldsymbol x, \kappa, \boldsymbol d_j)\|_{L^2(\partial B_R)}\lesssim \|u^s_i(\boldsymbol x, \kappa, \boldsymbol d_j)\|_{H^1(B_R)}\lesssim 1
\]
and
\[
\|\partial_\nu u^s_i(\boldsymbol x, \kappa, \boldsymbol d_j)\|_{L^2(\partial B_R)}\lesssim \|u^s_i(\boldsymbol x, \kappa, \boldsymbol d_j)\|_{H^2(B_R)}\lesssim \kappa.
\]
Consequently, we have for $j=1,2$
\[
\begin{split}
&\|u_j\|_{L^2(\partial B_R)},\,\,\|w_j\|_{L^2(\partial B_R)}=\mathcal{O}(1),\\
&\|\partial_\nu u_j\|_{L^2(\partial B_R)},\,\,\|\partial_\nu w_j\|_{L^2(\partial B_R)}=\mathcal{O}(\kappa).
\end{split}
\]

Denote
\[
\epsilon^2(\kappa,\boldsymbol d_1,\boldsymbol d_2) = I(\kappa,\boldsymbol d_1,\boldsymbol d_2)
\]
and
\begin{align*}
I(\kappa,\boldsymbol d_1,\boldsymbol d_2)&= \sum_{j = 1, 2}\kappa^2\|R_{V_1}(\kappa)(V_1 e^{\mathrm{i}\kappa \boldsymbol x\cdot \boldsymbol d_j})-R_{V_2}(\kappa)(V_2 e^{\mathrm{i}\kappa \boldsymbol x\cdot \boldsymbol d_j})\|^2_{L^2(\partial B_R)}\\
&\quad + \|\partial_\nu\Big(R_{V_1}(\kappa)(V_1 e^{\mathrm{i}\kappa \boldsymbol x\cdot \boldsymbol d_j})-R_{V_2}(\kappa)(V_2 e^{\mathrm{i}\kappa \boldsymbol x\cdot \boldsymbol d_j})\Big)\|^2_{L^2(\partial B_R)}, \quad \kappa>0.
\end{align*}
Combining above identities and estimates, we have
\begin{equation}\label{est1}
\left|\int_{B_R}V(u_1u_2+w_1w_2)\mathrm{d}x\right|^2\lesssim \epsilon^2(\kappa,\boldsymbol d_1,\boldsymbol d_2).
\end{equation}
Moreover, we have 
\begin{align}\label{iden4}
&\int_{B_R}Vu_1u_2\mathrm{d}\boldsymbol x\notag\\
=&\int_{B_R} Ve^{\mathrm{i}\kappa \boldsymbol x\cdot (\boldsymbol d_1+\boldsymbol d_2)}-Ve^{\mathrm{i}\kappa \boldsymbol x\cdot \boldsymbol d_1}R_{V_2}(\kappa)(V_2e^{\mathrm{i}\kappa \boldsymbol x\cdot \boldsymbol d_2})\notag\\
&\quad -Ve^{\mathrm{i}\kappa \boldsymbol x\cdot \boldsymbol d_2}R_{V_1}(\kappa)(V_1e^{\mathrm{i}\kappa \boldsymbol x\cdot \boldsymbol d_1})+VR_{V_1}(\kappa)(V_1e^{\mathrm{i}\kappa \boldsymbol x\cdot \boldsymbol d_1})R_{V_2}(\kappa)(V_2e^{\mathrm{i}\kappa \boldsymbol x\cdot \boldsymbol d_2})\mathrm{d}\boldsymbol x\notag\\
=&\hat{V}(-\kappa(\boldsymbol d_1 + \boldsymbol d_2)) + \text{remainders},
\end{align}
and similarly,
\begin{align}\label{iden5}
&\int_{B_R}Vw_1w_2\mathrm{d}\boldsymbol x\notag\\
=&\int_{B_R} Ve^{\mathrm{i}\kappa \boldsymbol x\cdot (\boldsymbol d_1+\boldsymbol d_2)}-Ve^{\mathrm{i}\kappa \boldsymbol x\cdot \boldsymbol d_1}R_{V_1}(\kappa)(V_1e^{\mathrm{i}\kappa \boldsymbol x\cdot \boldsymbol d_2})\notag\\
&\quad-Ve^{\mathrm{i}\kappa \boldsymbol x\cdot \boldsymbol d_2}R_{V_2}(\kappa)(V_2e^{\mathrm{i}\kappa \boldsymbol x\cdot \boldsymbol d_1})+VR_{V_1}(\kappa)(V_1e^{\mathrm{i}\kappa \boldsymbol x\cdot \boldsymbol d_2})R_{V_2}(\kappa)(V_2e^{\mathrm{i}\kappa \boldsymbol x\cdot \boldsymbol d_1})\mathrm{d}\boldsymbol x\notag\\
=&\hat{V}(-\kappa(\boldsymbol d_1 + \boldsymbol d_2)) + \text{remainders}.
\end{align}

Using the resolvent estimate in Lemma \ref{bound}, we obtain
\[
\|R_{V_i}(V_i e^{\mathrm{i}\kappa \boldsymbol x\cdot \boldsymbol d_j})\|_{L^2(B_R)}\lesssim\frac{1}{\kappa}, \quad 1\leq i, j\leq 2,
\]
and thus the remainders in both \eqref{iden4} and \eqref{iden5} satisfy 
\[
| \text{remainders}|\lesssim\frac{1}{\kappa}.
\]
We add \eqref{iden4} and \eqref{iden5} and use the estimate \eqref{est1} to have
\[
\left|\hat{V}(-\kappa(\boldsymbol d_1 + \boldsymbol d_2))\right|^2\lesssim \epsilon^2(\kappa,\boldsymbol d_1,\boldsymbol d_2)+\frac{1}{\kappa^2}.
\]
Thus, noting that 
\[
\{-\kappa(\boldsymbol d_1 + \boldsymbol d_2), \boldsymbol d_1\in\mathbb S^2, \boldsymbol d_2\in\mathbb S^2\} = \{\boldsymbol \xi: |\boldsymbol \xi|\leq 2\kappa\},
\]
we have, for all $|\xi|\leq 2\kappa$, that
\begin{align}\label{crucial}
|\hat{V}(\boldsymbol \xi)|^2\lesssim  \sup_{\boldsymbol d_1\in\mathbb S^2,\boldsymbol d_2\in\mathbb{S}^2}\epsilon^2(\kappa,\boldsymbol d_1,\boldsymbol d_2) + \frac{1}{\kappa^2}.
\end{align}

Noting
\[
\overline{R_{V_i}(\kappa) (V_i e^{\mathrm{i}\kappa \boldsymbol x\cdot \boldsymbol d_j})} = R_{V_i}(-\kappa) (V_i e^{-\mathrm{i}\kappa \boldsymbol x\cdot \boldsymbol d_j}),
\]
we may extend the domain of $I(\kappa)$ from $\mathbb R^+$ to $\mathbb C$ by
\[
\begin{split}
&I(\kappa,\boldsymbol d_1,\boldsymbol d_2)=\\
&\kappa^2\int_{\partial B_R} (R_{V_1}(\kappa)(V_1 e^{\mathrm{i}\kappa \boldsymbol x\cdot \boldsymbol d_j})-R_{V_2}(\kappa)(V_2 e^{\mathrm{i}\kappa \boldsymbol x\cdot \boldsymbol d_j}))
(R_{V_1}(-\kappa)(V_1 e^{-\mathrm{i}\kappa \boldsymbol x\cdot \boldsymbol d_j})-R_{V_2}(-\kappa)(V_2 e^{-\mathrm{i}\kappa \boldsymbol x\cdot \boldsymbol d_j}))\mathrm{d}s\\
&+\int_{\partial B_R} \partial_\nu(R_{V_1}(\kappa)(V_1 e^{\mathrm{i}\kappa \boldsymbol x\cdot \boldsymbol d_j})-R_{V_2}(\kappa)(V_2 e^{\mathrm{i}\kappa \boldsymbol x\cdot \boldsymbol d_j}))
\partial_\nu(R_{V_1}(-\kappa)(V_1 e^{-\mathrm{i}\kappa \boldsymbol x\cdot \boldsymbol d_j})-R_{V_2}(-\kappa)(V_2 e^{-\mathrm{i}\kappa \boldsymbol x\cdot \boldsymbol d_j}))\mathrm{d}s.
\end{split}
\]
Since $R_{V_j}(\kappa)$ is meromorphic in $\mathbb C$, we know that $I(\cdot,\boldsymbol d_1,\boldsymbol d_2)$ is a meromorphic function on $\mathbb C$. By Lemma \ref{bound} we have that $I(\kappa)$ is holomorphic in a domain 
\begin{equation}\label{regionR}
\mathcal{R} = \{z\in \mathbb C: (K_0, +\infty)\times (-d, d) \},
\end{equation}
for some constants $K_0>0$ and $d>0$.\\

Suppose that we are given the following near-field multi-frequency data discrepancy in a finite interval $I = [K_0, K]$ with $K>K_0>C_0$ where $C_0$ is specified in Lemma \ref{bound}
\begin{align*}
\epsilon^2 = \sup_{\kappa\in I,\boldsymbol d_1\in\mathbb S^2, \boldsymbol d_2\in\mathbb S^2}I(\kappa,d_1,d_2):=&\sup_{\kappa\in I, \boldsymbol d_1\in\mathbb S^2, \boldsymbol d_2\in\mathbb S^2}\Big(\sum_{j = 1, 2}\Big( \kappa^2\|u^s_1(\boldsymbol x, \kappa, \boldsymbol d_j) - u^s_2(\boldsymbol x, \kappa, \boldsymbol d_j)\|_{L^2(\partial B_R)}^2\\
&\quad\quad\quad\quad+\|\partial_\nu (u^s_1(\boldsymbol x, \kappa, \boldsymbol d_j) - u^s_2(\boldsymbol x, \kappa, \boldsymbol d_j))\|_{L^2(\partial B_R)}^2\Big)\Big)\\
=&\sup_{\kappa\in I, \boldsymbol d\in\mathbb S^2}2\Big( \kappa^2\|u^s_1(\boldsymbol x, \kappa, \boldsymbol d) - u^s_2(\boldsymbol x, \kappa, \boldsymbol d)\|_{L^2(\partial B_R)}^2\\
&\quad\quad\quad\quad+\|\partial_\nu (u^s_1(\boldsymbol x, \kappa, \boldsymbol d) - u^s_2(\boldsymbol x, \kappa, \boldsymbol d))\|_{L^2(\partial B_R)}^2\Big).
\end{align*}
It is shown in \cite{LY} that there exists a Dirichlet-to-Neumann (DtN) operator $T$ such that 
\[
\partial_\nu u^s = T u^s \quad \text{on} \,\, \partial B_R,
\]
which is the transparent boundary condition for the scattering problem of the Schr\"odinger equation. Thus one has
\[
\epsilon^2=\sup_{\kappa\in I, \boldsymbol d\in\mathbb S^2}2\Big( \kappa^2\|u^s_1(\boldsymbol x, \kappa, \boldsymbol d) - u^s_2(\boldsymbol x, \kappa, \boldsymbol d)\|_{L^2(\partial B_R)}^2+\|T (u^s_1(\boldsymbol x, \kappa, \boldsymbol d) - u^s_2(\boldsymbol x, \kappa, \boldsymbol d))\|_{L^2(\partial B_R)}^2\Big).
\]

The following lemma  \cite[Lemma A.1]{LZZ} gives a link between the values of an analytic function for small and large arguments.
\begin{lemma}\label{ac}
Let $p(z)$ be analytic in the infinite rectangular slab $\mathcal{R}$ as in \eqref{regionR},
and continuous in $\overline{\mathcal{R}}$ satisfying
\begin{align*}
\begin{cases}
|p(z)|\leq \epsilon, &\quad z\in (K_0, K],\\
|p(z)|\leq M, &\quad z\in \mathcal{R},
\end{cases}
\end{align*}
where $K_0, K, \epsilon$ and $M$ are positive constants. Then there exists a function $\mu(z)$ with $z\in (K, +\infty)$ satisfying 
\begin{equation}\label{LB}
\mu(z) \geq \frac{64ad}{3\pi^2(a^2 + 4d^2)} e^{\frac{\pi}{2d}(\frac{a}{2} - z)},
\end{equation}
where $a = K - K_0$, such that
\begin{align}\label{p}
|p(z)|\leq M\epsilon^{\mu(z)}\quad \forall\, z\in (K, +\infty).
\end{align}
\end{lemma}

Now, recall that our data discrepancy condition implies that
for all $\boldsymbol d_1,\boldsymbol d_2\in\mathbb{S}^2$
\[
|I(z,\boldsymbol d_1,\boldsymbol d_2)|\leq\epsilon^2,\quad z\in(K_0, K],
\]
and then
\[
|z^{-2}I(z,\boldsymbol d_1,\boldsymbol d_2)|\lesssim\epsilon^2,\quad z\in(K_0, K].
\]
By Lemma \ref{bound}, we have for $|\Im z|<d$
\[
\begin{split}
&\|R_{V_i}(\pm z)(V_ie^{\pm\mathrm{i}z \boldsymbol x\cdot \boldsymbol d_j})\|_{H^1(B_R)}\lesssim 1,\\
&\|R_{V_i}(\pm z)(V_ie^{\pm\mathrm{i}z \boldsymbol x\cdot \boldsymbol d_j})\|_{H^2(B_R)}\lesssim |z|.
\end{split}
\]
Therefore
\[
\begin{split}
&\|R_{V_i}(\pm z)(V_ie^{\pm\mathrm{i}z \boldsymbol x\cdot \boldsymbol d_j})\|_{L^2(B_R)}\lesssim 1,\\
&\|\partial_\nu(R_{V_i}(\pm z)(V_ie^{\pm\mathrm{i}z \boldsymbol x\cdot \boldsymbol d_j}))\|_{L^2(B_R)}\lesssim |z|,
\end{split}
\]
and consequently
\[
|z^{-2}I(z,\boldsymbol d_1,\boldsymbol d_2)|\leq M,\quad z\in\mathcal{R},
\]
where $M$ is a positive constant. Applying the quantitative analytic continuation principle in Lemma \ref{ac} to $p(z)=z^{-2}I(z,\boldsymbol d_1,\boldsymbol d_2)$, we have
\[
|\kappa^{-2}I(\kappa,\boldsymbol d_1,\boldsymbol d_2)|\lesssim \epsilon^{2\mu(\kappa)},\quad \kappa>K,
\]
where 
\[
\mu(\kappa)=\frac{64ad}{3\pi^2(a^2+4d^2)}e^{\frac{\pi}{2 d}(\frac{a}{2}-\kappa)}, \quad a = K - K_0,
\]
which gives that
\[
\kappa^{-2}\epsilon^2(\kappa,\boldsymbol d_1,\boldsymbol d_2) \lesssim \epsilon^{2\mu(\kappa)}, \quad \kappa>K.
\]
Fixing a number $A$, we notice that for $\kappa\in (K,A]$,
\[
\mu(\kappa)\gtrsim ce^{-\sigma A},
\]
and then
\[
|\kappa^{-2}\epsilon^2(\kappa,\boldsymbol d_1,\boldsymbol d_2)|\lesssim \exp\{-ce^{-\sigma A}|\ln \epsilon|\}.
\]
%where $c$ and $\sigma$ are some positive constants,
%$A = K_1 - \frac{K - K_0}{2}$ and $K_1$ is some positive constant to be determined later.
Using the fact 
$e^{-x}\leq\frac{6!}{x^{6}}$ for $x>0$, we obtain
\[
|\epsilon^2(\kappa,\boldsymbol d_1,\boldsymbol d_2)|\lesssim e^{6\sigma A}|\ln \epsilon|^{-6},\quad \kappa\in (K,A],
\]
for some constant $\sigma>0$.

To derive the increasing stability estimates, we consider the following two situations.

\textbf{Case 1:} $K\leq \frac{1}{2\sigma}\ln|\ln \epsilon|$.
Taking $A=\frac{1}{2\sigma}\ln|\ln \epsilon|$, we have, in particular,
\[
|\epsilon^2(A,\boldsymbol d_1,\boldsymbol d_2)|\lesssim |\ln\epsilon|^{-3}.
\]
By \eqref{crucial}, we have
\[
|\hat{V}(\boldsymbol \xi)|^2\leq |\ln\epsilon|^{-3}+\frac{1}{A^2}
\]
for all $\boldsymbol \xi$ satisfying $|\boldsymbol \xi|\leq 2A$.
Then we have
\[
\int_{|\boldsymbol \xi|\leq 2A^{\frac{2}{3+2s}}}|\hat{V}(\boldsymbol \xi)|^2\mathrm{d}\boldsymbol \xi\lesssim (\ln|\ln\epsilon|)^{\frac{6}{3+2s}}|\ln\epsilon|^{-3}+(\ln|\ln\epsilon|)^{-\frac{4s}{3+2s}}\lesssim \frac{1}{(\ln|\ln\epsilon|)^{2\beta}}.
\]
Here we denote $\beta=\frac{2s}{3+2s}$.
By the assumption $\|V\|_{H^s}\leq Q$, we have
\[
\int_{|\boldsymbol\xi|> 2A^{\frac{2}{3+2s}}}|\hat{V}(\boldsymbol\xi)|^2\mathrm{d}\xi\lesssim\frac{1}{A^{\frac{4s}{3+2s}}}\leq (\ln|\ln \epsilon|)^{-2\beta}.
\]
Combing above estimates, we obtain
\[
\begin{split}
\|V\|_{L^2}^2=&\int_{|\boldsymbol\xi|\leq 2A^{\frac{2}{3+2s}}}|\hat{V}(\boldsymbol\xi)|^2\mathrm{d}\boldsymbol\xi+\int_{|\boldsymbol\xi|> 2A^{\frac{2}{3+2s}}}|\hat{V}(\boldsymbol\xi)|^2\mathrm{d}\xi\\
\lesssim &(\ln|\ln \epsilon|)^{-2\beta}\\
\lesssim &\frac{1}{K^\beta(\ln|\ln \epsilon|)^\beta}.
\end{split}
\]
\textbf{Case 2:} $K\geq \frac{1}{2\sigma}\ln|\ln \epsilon|$. First notice that
\[
\int_{|\boldsymbol\xi|\leq 2K^{\frac{2}{3+2s}}}|\hat{V}(\boldsymbol\xi)|^2\mathrm{d}\boldsymbol\xi\lesssim K^{\frac{6}{3+2s}}\epsilon^2+K^{-2\beta}.
\]
Also, we use the estimate
\[
\int_{|\boldsymbol\xi|> 2K^{\frac{2}{3+2s}}}|\hat{V}(\boldsymbol\xi)|^2\mathrm{d}\boldsymbol\xi\lesssim \frac{1}{K^{\frac{4s}{3+2s}}}\leq K^{-2\beta}.
\]
Now we have
\[
\begin{split}
\|V\|_{L^2}^2=&\int_{|\boldsymbol\xi|\leq 2K^{\frac{2}{3+2s}}}|\hat{V}(\boldsymbol\xi)|^2\mathrm{d}\xi+\int_{|\boldsymbol\xi|> 2K^{\frac{2}{3+2s}}}|\hat{V}(\boldsymbol\xi)|^2\mathrm{d}\xi\\
\lesssim &K^{\frac{6}{3+2s}}\epsilon^2+K^{-2\beta}\\
\lesssim &K^{\frac{6}{3+2s}}\epsilon^2+\frac{1}{K^\beta(\ln|\ln \epsilon|)^\beta}.
\end{split}
\]
Combining Case 1 and Case 2 we prove the following theorem.

\begin{theorem}\label{main}
Let $V_1, V_2\in \mathcal C_Q$. The following increasing stability estimate holds
\begin{align}\label{stability}
\|V_1 - V_2\|_{L^2(B_R)}^2
\lesssim K^\alpha\epsilon^2+\frac{1}{K^\beta(\ln|\ln \epsilon|)^\beta},
\end{align}
where $\alpha = \frac{6}{3+2s}$ and $\beta = \frac{2s}{3 + 2s}$.
\end{theorem}

The stability estimate \eqref{stability} contains a Lipschitz type data discrepancy and a high wavenumber tail of the
potential function. Moreover, the latter decreases as the upper bound $K$ of the wavenumber increases. 
which makes the problem have an almost Lipschitz
stability. The result reveals that the problem becomes more stable when higher
frequency data is used.

\subsubsection{Far-field pattern}\label{farr}

In this section we study the increasing stability estimate using multi-frequency far-field pattern.

Now we consider the Schwartz kernel for the free resolvent $R_0(\kappa)$ (that is when $V=0$),
\[
R_0(\kappa,\boldsymbol x,\boldsymbol y)=\frac{e^{\mathrm{i}\kappa|\boldsymbol x-\boldsymbol y|}}{4\pi|\boldsymbol x-\boldsymbol y|},
\]
which is the three-dimensional fundamental solution of the free Schr\"odinger equation.
For $|\boldsymbol x|$ large, we have the asymptotics
\[
|\boldsymbol x-\boldsymbol y|=|\boldsymbol x|-\frac{\boldsymbol x\cdot \boldsymbol y}{|\boldsymbol x|}+\mathcal{O}\left(\frac{|\boldsymbol y|^2}{|\boldsymbol x|}\right).
\]
Then for fixed $\boldsymbol y$ and large $|\boldsymbol x|$, we have
\begin{align}\label{asymptotic}
R_0(\kappa,\boldsymbol x,\boldsymbol y)=\frac{1}{4\pi}|\boldsymbol x|^{-1}e^{\mathrm{i}\kappa|\boldsymbol x|}e^{-\mathrm{i}\kappa \frac{\boldsymbol x\cdot \boldsymbol y}{|\boldsymbol x|}}\left(1+\mathcal{O}\left(\frac{1}{|\boldsymbol x|}\right)\right).
\end{align}
Then we have the integral equation 
\[
u^s=\frac{e^{\mathrm{i}\kappa|\boldsymbol x|}}{4\pi|\boldsymbol x|}\int_{\mathbb{R}^3}e^{-\mathrm{i}\kappa \frac{\boldsymbol x\cdot \boldsymbol y}{|\boldsymbol x|}}Vu\mathrm{d}\boldsymbol y+\mathcal{O}\left(\frac{1}{|\boldsymbol x|^2}\right).
\]
Denote $\boldsymbol\theta = \frac{\boldsymbol x}{|\boldsymbol x|}.$ The far-field pattern $\mathcal{A}_{\infty}$ is defined by
\begin{equation}\label{bform}
\begin{split}
\mathcal{A}_{\infty}(\kappa,\boldsymbol\theta,\boldsymbol d)=
\frac{1}{4\pi}\int_{\mathbb{R}^3}e^{-\mathrm{i}\kappa \langle \boldsymbol y,\boldsymbol\theta\rangle}Vu\mathrm{d}\boldsymbol y.
\end{split}
\end{equation}
Using $u = e^{{\rm i}\kappa\boldsymbol x\cdot\boldsymbol d} + u^s$ and the resolvent estimate $\|u^s\|_{L^2(B_R)}\lesssim \frac{1}{\kappa}$, we have
\[
\mathcal{A}_{\infty}(\kappa,\boldsymbol \theta,\boldsymbol d)=\hat{V}(\kappa(\boldsymbol \theta-\boldsymbol d))+\mathcal{O}\left(\frac{1}{\kappa}\right).
\]

Next, given the incident field $e^{{\rm i}\kappa \boldsymbol x\cdot \boldsymbol d}$, assume that $\mathcal{A}_{\infty, 1}$ and $\mathcal{A}_{\infty, 2}$ are the far-field patterns corresponding to 
two potential functions $V_1$ and $V_2$, respectively. Denote $V = V_1 - V_2$. By above discussions one has
\begin{align*}
|\hat{V}(\kappa(\boldsymbol \theta - \boldsymbol d))| \lesssim |\mathcal{A}_{\infty, 1}(\kappa,\boldsymbol\theta,\boldsymbol d) - \mathcal{A}_{\infty, 2}(\kappa,\boldsymbol\theta,\boldsymbol d)| + \frac{1}{\kappa} ,
\end{align*}
and consequently
\begin{align*}
|\hat{V}(\boldsymbol\xi)|^2 \lesssim \sup_{\boldsymbol\theta, \boldsymbol d\in\mathbb S^2} |\mathcal{A}_{\infty, 1}(\kappa,\boldsymbol\theta,\boldsymbol d) - \mathcal{A}_{\infty, 2}(\kappa,\boldsymbol\theta,\boldsymbol d)|^2 + \frac{1}{\kappa^2} 
\end{align*}
for all $|\boldsymbol\xi|\leq 2\kappa$. We also refer to \cite{ABFG} for a similar estimate.
Denote 
\[
\epsilon^2(\kappa,\boldsymbol\theta,\boldsymbol d) =  |\mathcal{A}_{\infty, 1}(\kappa,\boldsymbol\theta,\boldsymbol d) - \mathcal{A}_{\infty, 2}(\kappa,\boldsymbol\theta,\boldsymbol d)|^2.
\]
Notice that $\mathcal{A}_{\infty, i}(\cdot,\boldsymbol\theta,\boldsymbol d)$ given in \eqref{bform} has meromorphic continuation to $\mathbb{C}$ and
\[
\mathcal{A}_{\infty, i}(-\kappa,\boldsymbol\theta,\boldsymbol d) =\overline{\mathcal{A}_{\infty, i}(\kappa,\boldsymbol\theta,\boldsymbol d)},\quad \kappa\in\mathbb{R},\quad i=1,2.
\]
Therefore we can extend the definition $\epsilon^2(\kappa,\boldsymbol\theta,\boldsymbol d)$ for complex $\kappa$ as
\[
\epsilon^2(\kappa,\boldsymbol \theta,\boldsymbol d) =(\mathcal{A}_{\infty, 1}(\kappa,\boldsymbol \theta,\boldsymbol d) - \mathcal{A}_{\infty, 2}(\kappa,\boldsymbol \theta,\boldsymbol d))(\mathcal{A}_{\infty, 1}(-\kappa,\boldsymbol \theta,\boldsymbol d) 
- \mathcal{A}_{\infty, 2}(-\kappa,\boldsymbol\theta,\boldsymbol d)).
\]
By Lemma \ref{bound} we have that $\epsilon^2(\cdot,\boldsymbol\theta,\boldsymbol d)$ is holomorphic in $\mathcal{R}$. As in Section \ref{near}, assume that we are given the multi-frequency data discrepancy in $I = [K_0, K]$ with $K_0>C_0$
as follows
\begin{align*}
\epsilon^2 = \sup_{\kappa\in I,\boldsymbol\theta, \boldsymbol d\in \mathbb{S}^2} \epsilon^2(\kappa,\boldsymbol\theta,\boldsymbol d).
\end{align*}
Then, we have
\[
|\epsilon^2(z,\boldsymbol\theta,\boldsymbol d)|\leq \epsilon^2,\quad z\in (K_0,K),
\]
and, by Lemma \ref{bound},
\[
|\epsilon^2(z,\boldsymbol\theta,\boldsymbol d)|\leq M,\quad z\in\mathcal{R},
\]
for some constant $M>0$.

Thus, following the arguments in Section \ref{near} we obtain the same increasing stability estimate.
\begin{theorem}\label{main1}
Let $V_1, V_2\in \mathcal C_Q$. The following increasing stability estimate holds
\begin{align}\label{stability1}
\|V_1 - V_2\|_{L^2(B_R)}^2
\lesssim K^\alpha\epsilon^2+\frac{1}{K^\beta(\ln|\ln \epsilon|)^\beta},
\end{align}
where $\alpha = \frac{6}{3+2s}$ and $\beta = \frac{2s}{3 + 2s}$.
\end{theorem}

\subsection{Two-dimensional case}\label{two}

In this section we discuss the inverse scattering problem in $\mathbb R^2$. In fact, the results in Section \ref{three} can be extended to $\mathbb R^2$
in a straightforward way using the resolvent estimates in two dimensions.

Denote a symmetric sectorial region in $\mathbb C$ by 
\[
S = \{\lambda: \arg\lambda\in [-\frac{\pi}{4}, \frac{\pi}{4}] \cup [\frac{3\pi}{4}, \frac{5\pi}{4}], \,\lambda\neq 0 \}.
\]
The following theorem provides a resonance-free region and resolvent estimates in this region for $R_V(\lambda)$ in two dimensions.
The proof is put in Appendix for brevity of the main context of this paper.

\begin{theorem}\label{meromorphic}
Let $\rho\in C_0^\infty(\mathbb R^2)$ with $\rho = 1$ on $\text{supp}V$. The resolvent $\rho R_V(\lambda)\rho$ is a meromorphic family of operators in $S$. Moreover, $\rho R_V(\lambda)\rho$ is analytic for $\lambda\in S\cap \Omega_\delta$
with the following resolvent estimates
\begin{align}\label{bound_3}
\|\rho R_V(\lambda)\rho\|_{L^2(\mathbb R^2)\rightarrow H^j(\mathbb R^2)}\leq C|\lambda|^{-1/2}(1 + |\lambda|^2)^{j/2}
e^{L(\Im\lambda)_-},\quad j = 0, 1, 2, 
\end{align}
where $t_{-}:=\max\{-t,0\}$, $L>{\rm diam}({\rm supp}\rho): = \sup\{|x - y| : x, y \in {\rm supp}\rho\}$. Here $\Omega_\delta$ is defined as
\begin{align*}
\Omega_\delta:=\{\lambda: {\Im}\lambda\geq -\delta {\rm log}(1 + |\lambda|), |\lambda|\geq C_0\},
\end{align*}
where $C_0$ is a positive constant and $\delta$ satisfies $0<\delta<\frac{1}{4L}$.  In particular, there are only finitely many poles in the domain
\[
\{\lambda: {\Im}\lambda\geq -\delta {\rm log}(1 + |\lambda|), |\lambda|\leq C_0\}\cap S.
\]
\end{theorem}
We have the following corollary from Theorem \ref{meromorphic}.
\begin{corollary}\label{R}
There exists an infinite slab
\begin{equation}\label{regionR}
\mathcal{R} = \{z\in \mathbb C: (K_0, +\infty)\times (-d, d) \}
\end{equation}
such that $\{\mathcal{R}\cup -\mathcal{R}\}\subset S\cap \Omega_\delta$ for some constants $K_0>0$ and $d>0$.
\end{corollary}

We introduce the far-field pattern in $\mathbb R^2$.
From the asymptotic expansion of the fundamental solution $H^{(1)}_0(\kappa|\boldsymbol x - \boldsymbol y|)$ of the two-dimensional free Schr\"odinger equation in \cite{CK} as $|\boldsymbol x|\to\infty$, we have that
\begin{align*}
u(\boldsymbol x, \kappa, \boldsymbol d) = e^{{\rm i}\kappa \boldsymbol x\cdot \boldsymbol d} + \frac{1 + \rm i}{4\sqrt{\pi}} \frac{e^{{\rm i}\kappa |\boldsymbol x|}}{|\boldsymbol x|^{1/2}\sqrt{\kappa}} \Big\{\int_{\mathbb R^2} e^{-{\rm i}\kappa \boldsymbol y\cdot \hat{\boldsymbol x}} V u {\rm d}\boldsymbol y
+ \mathcal{O}\Big(\frac{1}{|\boldsymbol x|}\Big)\Big\}.
\end{align*}
The far-field pattern $\mathcal{A}_\infty(\kappa, \hat{\boldsymbol x}, \boldsymbol d)$ is defined as follows
\begin{align}\label{far}
\mathcal{A}_\infty(\kappa, \hat{\boldsymbol x}, \boldsymbol d) = \int_{\mathbb R^2} e^{-{\rm i}\kappa \boldsymbol y\cdot \hat{\boldsymbol x}} V u {\rm d}\boldsymbol y
\end{align}
where $\hat{\boldsymbol x} = \frac{\boldsymbol x}{|\boldsymbol x|}\in\mathbb S$ is the observation angle. 
Since $u(\boldsymbol x, \kappa, \boldsymbol d) = e^{{\rm i}\kappa \boldsymbol x\cdot \boldsymbol d}  + u^s(\boldsymbol x, \kappa, \boldsymbol d)$ and by Theorem \ref{meromorphic} for $\kappa\geq C_0$
\[
\|u^s(\boldsymbol x, \kappa, \boldsymbol d)\|_{L^2(B_R)} = \mathcal{O} \Big(\frac{1}{\kappa^{1/2}}\Big),
\]
plugging the above expression of $u$ to \eqref{far} yields
\begin{align}\label{important}
\int_{\mathbb R^2} V e^{{\rm i}\kappa \boldsymbol y \cdot (\boldsymbol d - \hat{\boldsymbol x})} {\rm d}\boldsymbol y 
= \mathcal{A}_\infty(\kappa, \hat{\boldsymbol x}, \boldsymbol d)  + \mathcal{O}\Big(\frac{1}{\kappa^{1/2}}\Big),
\end{align} 
which gives the following crucial estimate for all $|\boldsymbol \xi|\leq 2\kappa$ similar to \eqref{crucial} that
\begin{align}\label{crucial_1}
|\hat{V}(\boldsymbol\xi)|^2\lesssim  \sup_{\boldsymbol d, \boldsymbol\theta \in\mathbb{S}} |\mathcal{A}_\infty(\kappa, \boldsymbol \theta, \boldsymbol d)|^2 + \frac{1}{\kappa}.
\end{align}

Let $\kappa\in(0,+\infty)$ and $\boldsymbol d_1\in\mathbb S, \boldsymbol d_2 \in \mathbb S$.
As before, denote the scattered field corresponding to the incident field $e^{\mathrm{i}\kappa \boldsymbol x\cdot \boldsymbol d_j}$
and $V_i$ by 
\begin{align}\label{data_2d}
u^s_i(\boldsymbol x, \kappa, \boldsymbol d_j) = -R_{V_i}(\kappa)(V_i e^{\mathrm{i}\kappa \boldsymbol x\cdot \boldsymbol d_j}), \quad 1\leq i, j \leq 2,
\end{align}
Assume that we are given the following near-field multi-frequency Dirichlet data discrepancy in a finite interval $I = [K_0, K]$ with $K>K_0>C_0$ where $C_0$ is specified in Theorem \ref{meromorphic}
\begin{align*}
\epsilon^2 
=&\sup_{\kappa\in I, \boldsymbol d\in\mathbb S}2\Big( \kappa^2\|u^s_1(\boldsymbol x, \kappa, \boldsymbol d) - u^s_2(\boldsymbol x, \kappa, \boldsymbol d)\|_{L^2(\partial B_R)}^2\\
&\quad\quad\quad\quad+\|T (u^s_1(\boldsymbol x, \kappa, \boldsymbol d) - u^s_2(\boldsymbol x, \kappa, \boldsymbol d))\|_{L^2(\partial B_R)}^2\Big).
\end{align*}

We have the following stability estimate in two dimensions by near-field multi-frequency data.
Notice that we have different index $\alpha$ and $\beta$ in $\mathbb R^2$ since from the resolvent estimate in \eqref{bound_3}
now we have $\|u\|_{L^2(B_R)}\lesssim \frac{1}{\sqrt{\kappa}}$ in $\mathbb R^2$ instead of $\|u\|_{L^2(B_R)}\lesssim \frac{1}{\kappa}$ in $\mathbb R^3$.

\begin{theorem}\label{main_2d}
Let $V_1, V_2\in \mathcal C_Q$. The following increasing stability estimate holds
\begin{align*}
\|V_1 - V_2\|_{L^2(B_R)}^2
\lesssim K^\alpha\epsilon^2+\frac{1}{K^\beta(\ln|\ln \epsilon|)^\beta},
\end{align*}
where $\alpha = \frac{3}{2(3+2s)}$ and $\beta = \frac{s}{2(3 + 2s)}$.
\end{theorem}

If one uses the far-field pattern as data, we let 
\begin{align*}
\epsilon^2 = \sup_{\kappa\in I, \boldsymbol d, \boldsymbol\theta\in \mathbb{S}} |\mathcal{A}_{1, \infty}(\kappa,\boldsymbol \theta, \boldsymbol d), - \mathcal{A}_{2, \infty}(\kappa,\boldsymbol \theta, \boldsymbol d)|^2.
\end{align*}
 We have the following increasing stability estimate.
 \begin{theorem}\label{main_2d_1}
Let $V_1, V_2\in \mathcal C_Q$. The following increasing stability estimate holds
\begin{align*}
\|V_1 - V_2\|_{L^2(B_R)}^2
\lesssim K^\alpha\epsilon^2+\frac{1}{K^\beta(\ln|\ln \epsilon|)^\beta},
\end{align*}
where $\alpha = \frac{3}{2(3+2s)}$ and $\beta = \frac{s}{2(3 + 2s)}$.
\end{theorem}

\section{The magnetic Schr\"odinger equation}\label{magnetic}

As an application of the method proposed in Section \ref{acoustic} for the classical Schr\"odinger equation,
in this section we continue to investigate the increasing stability of determining both the compactly supported unknown magnetic and electric potentials of the magnetic Schr\"odinger equation in three dimensions. Let $e^{\mathrm{i}\kappa {\boldsymbol x}\cdot {\boldsymbol d}}$ be an incident plane wave with direction $\boldsymbol d\in\mathbb S^{2}$.
We consider the scattering problem modeled by the following  magnetic Schr\"odinger equation in $\mathbb R^3$
\begin{equation}\label{eqn3}
H u - \kappa^2 u = 0,
\end{equation}
where $H$ denotes the magnetic Schr\"odinger operator 
\[
H =  -(\nabla + {\rm i} {\boldsymbol b}(\boldsymbol x))^2 + V(\boldsymbol x),
\]
$u(\boldsymbol x,\kappa, \boldsymbol d) = e^{\mathrm{i}\kappa \boldsymbol x\cdot \boldsymbol d} + u^s(\boldsymbol x,\kappa,\boldsymbol d)$ is the total field, $u^s$ is the outgoing scattered field, $\kappa>0$ is the wavenumber. Here $ {\boldsymbol b}(\boldsymbol x)\in  W_{\rm comp}^{1, \infty}(\mathbb R^3)^3$ is the vector magnetic potential function and $V(\boldsymbol x)\in L_{\rm comp}^\infty(\mathbb R^3)$ is the electric potential function. 
Since for the inverse scattering problem, it was shown in \cite{ER} that changing $\boldsymbol b$ to $\tilde{\boldsymbol b} = \boldsymbol b + \nabla \varphi$ does not change the far-field pattern if $\varphi$ satisfied some exponentially decay conditions.
Thus, in this paper we shall assume that $\boldsymbol b$ satisfies $\nabla\cdot \boldsymbol b = 0$.
We assume that both potential functions are real-valued.
We also assume that the
supports of ${\boldsymbol b}(\boldsymbol x)$ and $V(\boldsymbol x)$ are compactly contained in $B_R$.
For each $\kappa>0$ the direct scattering problem is well-posed \cite{Vodev}.

By rewriting the operator $H$ as follows
\[
H = -\Delta - 2{\rm i} {\boldsymbol b} \cdot \nabla + \tilde{V}
\]
where $\tilde{V} = |{\boldsymbol b}|^2 + V$,
the equation \eqref{eqn3} becomes
\begin{equation}\label{eqn4}
-\Delta u - 2{\rm i} {\boldsymbol b} \cdot \nabla u + \tilde{V} u - \kappa^2 u = 0.
\end{equation}
We are interested in the inverse scattering problem of determining the magnetic potential ${\boldsymbol b}$ and electric potential $V$ from the near-field Dirichlet boundary measurements $u^s(\boldsymbol x, \kappa, \boldsymbol d)\vert_{\partial B_R}$, $\boldsymbol x\in\partial B_R$, $\boldsymbol d\in \mathbb S^2$ or the far-field pattern defined by \eqref{mf} corresponding to the wavenumber $\kappa$ given in a finite interval.

We study the stability estimate for the magnetic Schr\"odinger equation. The proof applies the method proposed in Section \ref{acoustic}
to the more complicated magnetic Schr\"odinger equation.
\subsection{Near-field data}\label{near_1}
In this section we study the stability estimate of determining the potentials $\boldsymbol b$ and $V$ in the magnetic Schr\"odinger equation \eqref{eqn4} by multi-frequency near-field data.
Let $s>0$ be an arbitrary positive constant.
We introduce a vector real-valued functional space for the vector magnetic potential
\begin{align*}
\mathcal D_Q= \{\boldsymbol b \in H^{s}(\mathbb R^3)^3\cap W^{1, \infty}(\mathbb R^3)^3:  \nabla\cdot \boldsymbol b = 0, \|\boldsymbol b\|_{H^{s}(\mathbb R^3)^3}\leq Q, \|\boldsymbol b\|_{L^{\infty}(\mathbb R^3)^3}\leq Q,
 ~ \boldsymbol b: \mathbb R^3 \rightarrow \mathbb R^3\}.
\end{align*}
 We further assume that the $W^{1, \infty}$ norm of $\boldsymbol b$ in $\mathcal D_Q$  is sufficiently small such that Lemma \ref{bound_1} holds.
Now let $\boldsymbol b_1$ and $\boldsymbol b_2$ be two magnetic potentials, and
$V_1$ and $V_2$ be two electric potentials. 
Let $\kappa\in(0,+\infty)$ and $\boldsymbol d_1\in\mathbb S^2, \boldsymbol d_2 \in \mathbb S^2$.

Denote the outgoing resolvent of the elliptic operator $H$ by 
\begin{equation}\label{resolvent_1}
R_{\boldsymbol b, V}(\lambda) = (H - \lambda^2)^{-1}.
\end{equation}
First we notice from \cite{Dyatlov} that the free resolvent $R_{0, 0}(\lambda): L^2(\mathbb R^3)\to L^2(\mathbb R^3)$ is well defined for $\Im\lambda>0$. 
Moreover,  one has that $R_{0, 0}(\lambda)$ is analytic for $\lambda\in\mathbb C$ as the following family of operators
\[
R_{0, 0}(\lambda): L_{\text{comp}}^2(\mathbb{R}^3)\rightarrow H^j_{\text{loc}}(\mathbb{R}^3),\quad j = 0, 1, 2.
\]
Here
\[
L_{\text{comp}}^2(\mathbb{R}^3):\{u\in L^2(\mathbb{R}^3),\exists R>0, |x|>R\Rightarrow |u(x)|=0\},
\]
\[
H^j_{\text{loc}}(\mathbb{R}^3):\{u\in\mathcal{D}'(\mathbb{R}^3),\forall \chi\in C_c^\infty(\mathbb{R}^3),\chi u\in H^j(\mathbb{R}^3)\}.
\]
Moreover, one has the following resolvent estimates
\begin{align}\label{fr}
\|R_{0,0}(\lambda)\|_{L^2_{\rm comp}(\mathbb R^3)\to H^j_{\rm loc}(\mathbb R^3)} \lesssim e^{T(\Im\lambda)_-}(|\lambda|^2 + 1)^{\frac{j-1}{2}}, \quad j = 0, 1, 2.
\end{align}
Thus, in the presence of only the electric potential $V$, for $\Im\lambda\gg 1$ one has that
\[
R_{0, V}(\lambda) = R_{0, 0}(\lambda) (I + VR_0(\lambda))^{-1}: L^2(\mathbb R^3)\to L^2(\mathbb R^3).
\]
The condition $\Im\lambda\gg 1$ guarantees the invertibility of $(I + VR_0(\lambda))^{-1}$ by the resolvent estimate \eqref{fr} and a Neumann series argument.
Then using the analytic Fredholm theorem one can meromorphically continue $R_{0, V}(\lambda)$ from $\Im\lambda\gg 1$ to the low-half complex plane as the following
family of operators
\[
R_{0, V}(\lambda): L_{\text{comp}}^2(\mathbb{R}^3)\rightarrow H^j_{\text{loc}}(\mathbb{R}^3),\quad j = 0, 1, 2.
\]
Furthermore, in the presence of both the electric potential $V$ and magnetic potential $\boldsymbol b$, one has the resolvent identity
\[
R_{0, V}(\lambda) = R_{0, 0}(\lambda) (I + (\boldsymbol b\cdot\nabla + V)R_0(\lambda))^{-1}: L^2(\mathbb R^3)\to L^2(\mathbb R^3), \quad \Im\lambda\gg 1,
\]
assuming that $\boldsymbol b$ is sufficiently small. In a similar way one can meromorphically continue $R_{\boldsymbol b, V}(\lambda)$ from $\Im\lambda\gg 1$ to the low-half complex plane as the following family of
operators
\[
R_{\boldsymbol b, V}(\lambda): L_{\text{comp}}^2(\mathbb{R}^3)\rightarrow H^j_{\text{loc}}(\mathbb{R}^3),\quad j = 0, 1, 2.
\]
%Based on the perturbation arguments in \cite{Dyatlov},
%we have that $R_{\boldsymbol b, V}(\lambda)$ is meromorphic in the complex plane
%for $V\in L^\infty_{\rm comp}(\mathbb R^3)$ and small magnetic potential $\boldsymbol b \in W^{1, \infty}(\mathbb R^3)^3$ as the following family of operators
%\[
%R_{\boldsymbol b, V}(\lambda): L_{\text{comp}}^2(\mathbb{R}^3)\rightarrow L^2_{\text{loc}}(\mathbb{R}^3),\quad \lambda\in\mathbb{C}.
%\]
%As indicated in \cite{Salo, Vodev}, for small magnetic potential $\boldsymbol b \in W^{1, \infty}(\mathbb R^3)^3$, 
Moreover, an analytic domain and
resolvent estimates as those in Lemma \ref{bound} can be derived for the resolvent $R_{\boldsymbol b, V}(\lambda)$.
For details on the above arguments we refer the reader to \cite[Theorem 3.8 and 3.10]{Dyatlov} and their proofs.
In summary, we have the following lemma.
\begin{lemma}\label{bound_1}
Fix a cutoff function $\rho\in C_0^\infty(\mathbb R^3)$. Suppose $\boldsymbol b\in W^{1, \infty}(\mathbb R^3)^3$ with $\|\boldsymbol b\|_{W^{1, \infty}(\mathbb R^3)^3}$ sufficiently small and $V\in L^\infty(\mathbb{R}^3)$ with compact supports.
There exist constants $A$, $C, C_0$ and $T$ such that, for $j=0,1,2$,
\begin{align}
\|\rho R_{\boldsymbol b, V}(\lambda)\rho\|_{L^2(\mathbb R^3) \to H^j(\mathbb R^3)}\leq
C|\lambda|^{j - 1} e^{T(\mathrm{Im}\lambda)_-},
\end{align}
for 
\[
\Im\lambda\geq -A-\delta\log(1+|\lambda|),\quad|\lambda|>C_0,
\]
where $\delta<1/\mathrm{diam}(\mathrm{supp} V)$.

\end{lemma}

Since the direct scattering problem is well-posed for $\kappa>0$ (cf \cite{ER}), 
one has that $R_{\boldsymbol b, V}(\kappa)$ is well defined for $\kappa>0$.
Then from equation \eqref{eqn4} we can denote the scattered field corresponding to the incident field $e^{\mathrm{i}\kappa x\cdot \boldsymbol d_j}$, $\boldsymbol b_i$
and $V_i$ by 
\begin{align}\label{data}
u^s_i(\boldsymbol x, \kappa, \boldsymbol d_j) = -R_{\boldsymbol b_i, V_i}(\kappa)(2{\rm i}\boldsymbol b_i\cdot(\nabla e^{\mathrm{i}\kappa \boldsymbol x\cdot \boldsymbol d_j}) + \tilde{V}_i e^{\mathrm{i}\kappa \boldsymbol x\cdot \boldsymbol d_j}), \quad 1\leq i, j \leq 2,
\end{align}
where $\tilde{V}_i = |\boldsymbol b_i|^2 + V_i $.

Denote
\begin{align*}
u_j&=e^{\mathrm{i}\kappa \boldsymbol x\cdot \boldsymbol d_j}-R_{\boldsymbol b_j, V_j}(\kappa)(2{\rm i}\boldsymbol b_j\cdot(\nabla e^{\mathrm{i}\kappa \boldsymbol x\cdot \boldsymbol d_j}) + \tilde{V}_j e^{\mathrm{i}\kappa \boldsymbol x\cdot \boldsymbol d_j}), \quad j = 1, 2,
%&= e^{\mathrm{i}\kappa x\cdot d_j} + 2\kappa R_{b_j, V_j}(\kappa) (b_j \cdot d_j e^{\mathrm{i}\kappa x\cdot d_j}) - R_{b_j, V_j}(\kappa) (\tilde{V}_j e^{\mathrm{i}\kappa x\cdot d_j}) ,\quad j=1,2,
\end{align*}
which satisfy
\begin{align}\label{eqn5}
H u_j -\kappa^2 u_j = 0, \quad j = 1, 2.
\end{align}
Let $\boldsymbol b = \boldsymbol b_1 - \boldsymbol b_2$ and $V = V_2 - V_1$. Subtracting \eqref{eqn5} when $j = 1$ by itself when $j = 2$ gives
\begin{align}\label{id1}
-\Delta u - 2{\rm i} \boldsymbol b_1 \cdot \nabla u + \tilde{V}_1 u - \kappa^2 u = -2{\rm i} \boldsymbol b \cdot\nabla u_2 + (|\boldsymbol b_2|^2 - |\boldsymbol b_1|^2 + V) u_2.
\end{align}
Multiplying both sides of \eqref{id1} by $u_1$, noting $\nabla\cdot \boldsymbol b_1 = 0$ and integrating by parts yields
\begin{align*}
\int_{B_R} -2{\rm i} \boldsymbol b \cdot\nabla u_2 u_1 + (|\boldsymbol b_2|^2 - |\boldsymbol b_1|^2 + V - {\rm i}\nabla\cdot \boldsymbol b) u_2 u_1 = \int_{\partial B_R} u_1 \partial_\nu  u_2 - u_2 \partial_\nu  u_1.
\end{align*}
Let 
\[
\begin{split}
w_1&=e^{\mathrm{i}\kappa \boldsymbol x\cdot \boldsymbol d_2}-R_{\boldsymbol b_1, V_1}(\kappa)(2{\rm i}\boldsymbol b_1\cdot(\nabla e^{\mathrm{i}\kappa \boldsymbol x\cdot \boldsymbol d_2}) + \tilde{V}_1 e^{\mathrm{i}\kappa \boldsymbol x\cdot \boldsymbol d_2}),\\
w_2&=e^{\mathrm{i}\kappa \boldsymbol x\cdot \boldsymbol d_1}-R_{\boldsymbol b_2, V_2}(\kappa)(2{\rm i}\boldsymbol b_2\cdot(\nabla e^{\mathrm{i}\kappa \boldsymbol x\cdot \boldsymbol d_1}) + \tilde{V}_2 e^{\mathrm{i}\kappa \boldsymbol x\cdot \boldsymbol d_1}),
\end{split}
\]
and $w = w_1 - w_2$. Similarly one has
\begin{align*}
\int_{B_R} -2{\rm i} \boldsymbol b \cdot\nabla w_2 w_1 + (|\boldsymbol b_2|^2 - |\boldsymbol b_1|^2 + V) w_2 w_1 = \int_{\partial B_R} w_1 \partial_\nu  w_2 - w_2 \partial_\nu  w_1.
\end{align*}
Following the arguments in Section \ref{acoustic} one has
\begin{align}\label{identity}
&\int_{B_R} -2{\rm i} \boldsymbol b \cdot\nabla u_2 u_1 + (|\boldsymbol b_2|^2 - |\boldsymbol b_1|^2 + V )u_1 u_2 \notag\\
&\quad + \int_{B_R} -2{\rm i} \boldsymbol b \cdot\nabla w_2 w_1 + (|\boldsymbol b_2|^2 - |\boldsymbol b_1|^2 + V )w_1 w_2\notag\\
&= \int_{\partial B_R} u_1 \partial_\nu  u_2 - u_2 \partial_\nu  u_1 + \int_{\partial B_R} w_1 \partial_\nu  w_2 - w_2 \partial_\nu  w_1,
\end{align}
where the boundary integrals above have the same form as \eqref{identity1} and thus can be rewritten as a sum of \eqref{iden1}--\eqref{iden3}.

Now we consider the first integral in \eqref{identity}
\[
\int_{B_R} -2{\rm i} \boldsymbol b \cdot\nabla u_2 u_1{\rm d}\boldsymbol x.
\]
As $\nabla\cdot \boldsymbol b = 0$, by taking Fourier transform one has $\boldsymbol \xi\cdot \hat{\boldsymbol b}(\boldsymbol\xi) = 0$ which gives
\[
\int_{B_R} \boldsymbol b \cdot (\boldsymbol d_1 + \boldsymbol d_2) e^{\mathrm{i}\kappa \boldsymbol x\cdot (\boldsymbol d_1 + \boldsymbol d_2)} = 0.
\]
Thus, we have
\begin{align*}
\int_{B_R} -2{\rm i} \boldsymbol b \cdot\nabla u_2 u_1 &= \int_{B_R} -2{\rm i} \boldsymbol b \cdot (\nabla e^{\mathrm{i}\kappa \boldsymbol x\cdot \boldsymbol d_2}) e^{\mathrm{i}\kappa \boldsymbol x\cdot \boldsymbol d_1} + \text{remainders}\\
&= 2\kappa  \int_{B_R} \boldsymbol b\cdot \boldsymbol d_2 e^{\mathrm{i}\kappa \boldsymbol x\cdot (\boldsymbol d_1 + \boldsymbol d_2)} + \text{remainders} \\
&= \kappa\int_{B_R} \boldsymbol b \cdot (\boldsymbol d_2 - \boldsymbol d_1) e^{\mathrm{i}\kappa \boldsymbol x\cdot (\boldsymbol d_1 + \boldsymbol d_2)} + \text{remainders},
\end{align*}
where we have used
\[
-\int_{B_R} \boldsymbol b \cdot \boldsymbol d_1 e^{\mathrm{i}\kappa \boldsymbol x\cdot (\boldsymbol d_1 + \boldsymbol d_2)} = \int_{B_R} \boldsymbol b \cdot \boldsymbol d_2 e^{\mathrm{i}\kappa \boldsymbol x\cdot (\boldsymbol d_1 + \boldsymbol d_2)}
\]
in the last identity. Using the resolvent estimate in Lemma \ref{bound_1} one has that
\[
\text{remainders} = \mathcal{O}(1).
\]

Denote the unit vector
\[
\boldsymbol v = \frac{\boldsymbol d_2 - \boldsymbol d_1}{|\boldsymbol d_2 - \boldsymbol d_1|}.
\]
Denote the angle between two vectors $\boldsymbol d_1$ and $\boldsymbol d_2$ by $\langle \boldsymbol d_1, \boldsymbol d_2\rangle$.
Next we show that
\begin{align}\label{span}
\Big\{\Big(\int_{B_R} \boldsymbol b \cdot \boldsymbol v \, e^{\mathrm{i}\kappa \boldsymbol x\cdot (\boldsymbol d_1 + \boldsymbol d_2)}{\rm d}\boldsymbol x\Big) \boldsymbol v: \boldsymbol d_1\in\mathbb S^2, \boldsymbol d_2 \in \mathbb S^2,
\langle \boldsymbol d_1, \boldsymbol d_2\rangle\geq\frac{\pi}{2}\Big\} 
= \{\hat{\boldsymbol b}(\boldsymbol \xi): |\boldsymbol \xi|\leq \sqrt{2}\kappa\}.
\end{align}
The assumption $\langle \boldsymbol d_1, \boldsymbol d_2\rangle\geq\frac{\pi}{2}$ ensures that $|\boldsymbol d_1 - \boldsymbol d_2|$ has a lower bound with $|\boldsymbol d_1 - \boldsymbol d_2|\geq \sqrt{2}$. To this end,
fix $\boldsymbol \xi$ with $|\boldsymbol \xi|\leq\sqrt{2}\kappa$.
We can choose two suitable pairs of $\boldsymbol d_1, \boldsymbol d_2$ and $\tilde{\boldsymbol d}_1, \tilde{\boldsymbol d}_2$ such that $\boldsymbol \xi = \boldsymbol d_1 + \boldsymbol d_2 = \tilde{\boldsymbol d}_1 + \tilde{\boldsymbol d}_2$
and the unit vectors $\boldsymbol e = \frac{\boldsymbol d_1 + \boldsymbol d_2}{|\boldsymbol d_1 + \boldsymbol d_2|}$, $\boldsymbol v = \frac{\boldsymbol d_2 - \boldsymbol d_1}{|\boldsymbol d_2 - \boldsymbol d_1|}$ and 
$\tilde{\boldsymbol v} = \frac{\tilde{\boldsymbol d}_2 - \tilde{\boldsymbol d}_1}{|\tilde{\boldsymbol d}_2 - \tilde{\boldsymbol d}_1|}$ are mutually orthogonal.
As a consequence, one has
\[
\int_{B_R} \boldsymbol b \cdot \frac{\boldsymbol d_2 - \boldsymbol d_1}{|\boldsymbol d_2 - \boldsymbol d_1|} e^{\mathrm{i}\kappa \boldsymbol x\cdot (\boldsymbol d_1 + \boldsymbol d_2)} = \hat{\boldsymbol b}(\boldsymbol \xi)\cdot \boldsymbol v
\]
and
\[
\int_{B_R} \boldsymbol b \cdot \frac{\tilde{\boldsymbol d}_2 - \tilde{\boldsymbol d}_1}{|\tilde{\boldsymbol d}_2 - \tilde{\boldsymbol d}_1|} e^{\mathrm{i}\kappa \boldsymbol x\cdot (\boldsymbol d_1 + \boldsymbol d_2)} = \hat{\boldsymbol b}(\boldsymbol \xi)\cdot \tilde{\boldsymbol v}.
\]
Moreover, as $\nabla \cdot\boldsymbol b = 0$ gives $\boldsymbol e\cdot\hat{\boldsymbol b} = 0$ one has
\[
\hat{\boldsymbol b} = (\hat{\boldsymbol b}\cdot \boldsymbol v) \boldsymbol v + (\hat{\boldsymbol b}\cdot \tilde{\boldsymbol v})  \tilde{\boldsymbol v},
\]
which proves \eqref{span}.

Let $K_0>C_0$ where $C_0$ is specified in Lemma \ref{bound_1}. Denote $I = [K_0, K]$ and
\begin{align*}
\epsilon^2 = \sup_{\kappa\in I,\boldsymbol d_1\in\mathbb S^2, \boldsymbol d_2\in\mathbb S^2}  I(\kappa,\boldsymbol d_1,\boldsymbol d_2):=&\sup_{\kappa\in I, \boldsymbol d_1\in\mathbb S^2, \boldsymbol d_2\in\mathbb S^2}
\Big(\sum_{j = 1, 2}\Big( \kappa^2\|u^s_1(\boldsymbol x, \kappa, \boldsymbol d_j) - u^s_2(\boldsymbol x, \kappa, \boldsymbol d_j)\|_{L^2(\partial B_R)}^2\\
&\quad\quad\quad\quad+\|T (u^s_1(\boldsymbol x, \kappa, \boldsymbol d_j) - u^s_2(\boldsymbol x, \kappa, \boldsymbol d_j))\|_{L^2(\partial B_R)}^2\Big)\Big),
\end{align*}
where $T$ is the transparent boundary condition,
and 
\[
\tilde{\epsilon}^2 = \sup_{\kappa\in I,\boldsymbol d_1\in\mathbb S^2, \boldsymbol d_2\in\mathbb S^2} \frac{1}{\kappa^2}  I(\kappa,\boldsymbol d_1,\boldsymbol d_2).
\]
Using the resolvent estimates in Lemma \ref{bound_1} the identity \eqref{identity} becomes 
\begin{align}\label{identity3}
2\kappa\int_{B_R} \boldsymbol b \cdot (\boldsymbol d_2 - \boldsymbol d_1) e^{\mathrm{i}\kappa \boldsymbol x\cdot (\boldsymbol d_1 + \boldsymbol d_2)} + \mathcal{O}(1) = \text{boundary integrals on} \, \partial B_R.
\end{align}
Then combing \eqref{span} and \eqref{identity3} we arrive at
for all $|\boldsymbol\xi|\leq \sqrt{2}\kappa$, that
\begin{align}\label{crucial_1}
|\hat{\boldsymbol b}(\boldsymbol \xi)|^2&\lesssim  \sup_{\boldsymbol d_1\in\mathbb S^2, \boldsymbol d_2\in\mathbb{S}^2}\frac{1}{\kappa^2}\epsilon^2(\kappa,\boldsymbol d_1,\boldsymbol d_2) + \frac{1}{\kappa^2}\notag\\
&\lesssim \sup_{\boldsymbol d_1\in\mathbb S^2, \boldsymbol d_2\in\mathbb{S}^2}\tilde{\epsilon}^2(\kappa,\boldsymbol d_1,\boldsymbol d_2) + \frac{1}{\kappa^2}.
\end{align}
Then following the arguments in Section \ref{near} we obtain an increasing stability estimate for the magnetic potential $\boldsymbol b$ of the same form as \eqref{stability} in Theorem \ref{main}.
Moreover, once we recover $\boldsymbol b$, the stability estimate for $V$ follows from the arguments in Section \ref{near} as well. In summary, we have the following stability estimates for $\boldsymbol b$ and $V$.
\begin{theorem}\label{main_1}
Let $\boldsymbol b_1, \boldsymbol b_2\in \mathcal D_Q$ and $V_1, V_2\in \mathcal C_Q$. The following increasing stability estimate holds
\begin{align*}
\|\boldsymbol b_1 - \boldsymbol b_2\|_{L^2(B_R)^3}^2
&\lesssim K^\alpha\tilde{\epsilon}^2+\frac{1}{K^\beta(\ln|\ln \tilde{\epsilon}|)^\beta},\\
 \|V_1 - V_2\|_{L^2(B_R)}^2
 &\lesssim K^\alpha\epsilon^2+\frac{1}{K^\beta(\ln|\ln \epsilon|)^\beta},
\end{align*}
where $\alpha = \frac{6}{3+2s}$ and $\beta = \frac{2s}{3 + 2s}$.
\end{theorem}

%\begin{remark}
%The proof of Theorem \ref{main_1} can be extended to the two-dimensional case once the resolvent estimates in \cite[Theorem 1.1]{Vodev} also hold in $\mathbb R^2$.
%\end{remark}

\subsection{Far-field pattern}

In this section we study the stability estimate of determining the potentials $\boldsymbol b$ and $V$ using multi-frequency far-field measurements.
We introduce the far-field pattern for the magnetic Schr\"odinger equation. The fundamental solution for the free Schr\"odinger equation has the form
\[
R_0(\kappa, \boldsymbol x, \boldsymbol y) = \frac{e^{\mathrm{i}\kappa|\boldsymbol x-\boldsymbol y|}}{4\pi|\boldsymbol x-\boldsymbol y|}
\] 
which satisfies 
\[
-\Delta R_0(\kappa, \boldsymbol x,\boldsymbol y) - \kappa^2  R_0(\kappa, \boldsymbol x, \boldsymbol y) = \delta(\boldsymbol x - \boldsymbol y).
\]
Then we have the integral equation
\begin{equation*}
u^s = \int_{B_R} R_0 (\kappa, \boldsymbol x, \boldsymbol y) (2{\rm i} \boldsymbol b\cdot\nabla u - \tilde{V}u) {\rm d}\boldsymbol y.
\end{equation*}
By the asymptotic behavior \eqref{asymptotic} of $R_0 (\kappa, \boldsymbol x, \boldsymbol y)$ we have
\begin{equation*}
u^s =  \frac{e^{{\rm i}\kappa|\boldsymbol x|}}{4\pi|\boldsymbol x|}\int_{B_R} e^{-{\rm i}\kappa\boldsymbol\theta\cdot \boldsymbol y} (2{\rm i} \boldsymbol b\cdot\nabla (e^{{\rm i}\kappa \boldsymbol y\cdot \boldsymbol d})  - \tilde{V} e^{{\rm i}\kappa \boldsymbol y\cdot \boldsymbol d} ){\rm d}\boldsymbol y + \mathcal{O}\Big(\frac{1}{|\boldsymbol x|^2}\Big),
\end{equation*}
where $\theta = \frac{\boldsymbol x}{|\boldsymbol x|}$. Following \cite{CK} we define the far-field pattern by
\begin{equation}\label{mf}
\mathcal{A}_\infty (\kappa, \boldsymbol\theta, \boldsymbol d) = \int_{B_R} e^{-{\rm i}\kappa\boldsymbol\theta\cdot \boldsymbol y} \Big(2{\rm i} \boldsymbol b\cdot\nabla (e^{{\rm i}\kappa \boldsymbol y\cdot \boldsymbol d})  - \tilde{V} e^{{\rm i}\kappa \boldsymbol y\cdot \boldsymbol d}\Big){\rm d}\boldsymbol y. 
\end{equation}
Using the resolvent estimate in Lemma \ref{bound_1} we have
\begin{equation*}
\mathcal{A}_\infty (\kappa, \boldsymbol\theta, \boldsymbol d) = \int_{B_R} -2\kappa e^{-{\rm i}\kappa\boldsymbol\theta\cdot \boldsymbol y} \boldsymbol b\cdot \boldsymbol d  e^{{\rm i}\kappa \boldsymbol y\cdot \boldsymbol d} {\rm d}\boldsymbol y + \mathcal{O}(1).
\end{equation*}

Given the incident field $e^{{\rm i}\kappa \boldsymbol x\cdot \boldsymbol d}$, assume that $\mathcal{A}_{\infty, 1}$ and $\mathcal{A}_{\infty, 2}$ are the far-field patterns corresponding to two pairs of potential functions $\boldsymbol b_1, V_1$ and $\boldsymbol b_2, V_2$, respectively.
Denote 
\[
\epsilon^2(\kappa,\boldsymbol\theta,\boldsymbol d) =  |\mathcal{A}_{\infty, 1}(\kappa,\boldsymbol\theta,\boldsymbol d) - \mathcal{A}_{\infty, 2}(\kappa,\boldsymbol\theta,\boldsymbol d)|^2
\]
and the multi-frequency data discrepancy for $\kappa\in I = [K_0, K]$ with $K_0>C_0$ by
\begin{align*}
\epsilon^2 = \sup_{\kappa\in I,\boldsymbol\theta\in\mathbb S^2, \boldsymbol d\in \mathbb{S}^2} \epsilon^2(\kappa,\boldsymbol\theta,\boldsymbol d).
\end{align*}
and 
\[
\tilde{\epsilon}^2 = \sup_{\kappa\in I,\boldsymbol\theta\in\mathbb S^2, \boldsymbol d\in \mathbb{S}^2}\frac{1}{\kappa^2} \epsilon^2(\kappa,\boldsymbol \theta,\boldsymbol d).
\]

Following the arguments in Section \ref{near_1} and \ref{near} we arrive at the following increasing stability estimate.
\begin{theorem}\label{main_2}
Let $\boldsymbol b_1, \boldsymbol b_2\in\mathcal D_Q$ and $V_1, V_2\in \mathcal C_Q$. The following increasing stability estimate holds
\begin{align*}
\|\boldsymbol b_1 - \boldsymbol b_2\|_{L^2(B_R)^3}^2
&\lesssim K^\alpha\tilde{\epsilon}^2+\frac{1}{K^\beta(\ln|\ln \tilde{\epsilon}|)^\beta},\\
 \|V_1 - V_2\|_{L^2(B_R)}^2
 &\lesssim K^\alpha\epsilon^2+\frac{1}{K^\beta(\ln|\ln \epsilon|)^\beta},
\end{align*}
where $\alpha = \frac{6}{3+2s}$ and $\beta = \frac{2s}{3 + 2s}$.
\end{theorem}

\section{Conclusion}

In this paper we present a unified stability theory of the inverse potential scattering problems for both the classical and magnetic Schr\"odinger
equations. 
For the classical Schr\"odinger equation in two and three dimensions, the increasing stability is achieved to reconstruct the potential function. For the magnetic Schr\"odinger equation, the increasing stability is obtained to reconstruct both the magnetic and electric potentials.
The stability estimates consist of the data discrepancy and the high frequency tail for the unknown potential function. The result shows that
the ill-posedness of the inverse potential scattering problems decreases as the frequency increases for the data.
The analysis requires either near-field data or far-field pattern at multiple frequencies. As multi-frequency data is available, the method does not resort to the construction of CGO solutions and can be utilized to study the two-dimensional case.

\appendix

\section{Proof of Theorem \ref{meromorphic}}

In this section we consider the meromorphic continuation of the following outgoing  resolvent in two dimensions
\[
R_V(\lambda) = (-\Delta - \lambda^2 + V)^{-1}.
\]
We derive a resonance-free region and resolvent estimates in this region for $R_V(\lambda)$.

We start with the free resolvent $R_0(\lambda)$ for $V\equiv 0$.
One has that the kernel of the outgoing free resolvent $R_0(\lambda) = (-\Delta - \lambda^2)^{-1}$
is the Hankel function of the first kind $H_0^{(1)}(\lambda|x - y|)$ in two dimensions.
Now we analyze the kernel $H_0^{(1)}(\lambda|x - y|)$.
From \cite{FY06} the Hankel function $H_0^{(1)}(\lambda|x - y|)$ has the following integral form
\begin{align}\label{kernel}
R_0 (\lambda,   x,     y) =  C e^{{\rm i} \lambda | x -  y|} \int_0^\infty e^{-t} t^{-\frac{1}{2}} \Big( \frac{t}{2} - {\rm i}
\lambda |    x -     y| \Big)^{-\frac{1}{2}} {\rm d}t,
\end{align}
where $C$ is a positive constant. 
In what follows, we show that $R_0(\lambda)$ is analytic for $\lambda\in S$
as the following operator
\[
R_V(\lambda): L_{\text{comp}}^2(\mathbb{R}^2)\rightarrow L^2_{\text{loc}}(\mathbb{R}^2),\quad \lambda\in S,
\]
where
\[
L_{\text{comp}}^2(\mathbb{R}^2):\{u\in L^2(\mathbb{R}^2),\exists R>0, |x|>R\Rightarrow |u(x)|=0\},
\]
\[
L^2_{\text{loc}}(\mathbb{R}^2):\{u\in\mathcal{D}'(\mathbb{R}^2),\forall \chi\in C_c^\infty(\mathbb{R}^3),\chi u\in L^2(\mathbb{R}^2)\}.
\]
Here $S$ is a symmetric sectorial region in $\mathbb C$ denoted by 
\[
S = \{\lambda: \arg\lambda\in [-\frac{\pi}{4}, \frac{\pi}{4}] \cup [\frac{3\pi}{4}, \frac{5\pi}{4}], \,\lambda\neq 0 \}.
\]
For $\lambda\in S$ a simple calculation yields 
\begin{align}\label{est}
|R_0 (\lambda,     x,     y)| \lesssim  \frac{e^{-\Im\lambda |  x -   y|}}{|\Re\lambda|^{\frac{1}{2}}|    x -     y|^{\frac{1}{2}}} \int_0^\infty e^{-t} t^{-\frac{1}{2}} {\rm d}t \lesssim  
\frac{e^{-\Im\lambda|    x -     y|}}{|\lambda|^{\frac{1}{2}}|    x -     y|^{\frac{1}{2}}}.
\end{align}

The following theorem concerns the analytic continuation of the free resolvent $R_0(\lambda)$. 

\begin{theorem}\label{free_estimate_2d}
The free resolvent $R_0(\lambda)$ is analytic for $\lambda\in S, \Im\lambda>0$ as a family of operators
\begin{align*}
R_0 (\lambda): L^2(\mathbb R^{2})\to L^2(\mathbb R^{2})
\end{align*}
where $\|R_0 (\lambda)\|_{L^2(\mathbb R^{2})\to L^2(\mathbb R^{2})} = \mathcal{O}(1/|\lambda|^{1/2})$.
Moreover, for each $\rho\in C_0^\infty(\mathbb R^{2})$ the free resolvent $R_0(\lambda)$ extends to a family of analytic operators for $\lambda\in S$ as follows
\begin{align*}
\rho R_0 (\lambda) \rho: L^2(\mathbb R^{2})\to L^2(\mathbb R^{2})
\end{align*}
with the resolvent estimates
\begin{align}\label{free}
\|\rho R_0(\lambda) \rho\|_{L^2(\mathbb R^2)\rightarrow H^j(\mathbb R^2)}\lesssim |\lambda|^{-\frac{1}{2}} (1+|\lambda|^2)^{\frac{j}{2}} e^{L (\Im\lambda)_-}, \quad j=0, 1, 2,
\end{align}
where $t_{-}:=\max\{-t,0\}$ and $L>{\rm diam}({\rm supp}\rho): = \sup\{|  x -   y| :   x,   y \in {\rm supp}\rho\}$.
\end{theorem}

\begin{proof}

We first prove that $\rho R_0(\lambda)\rho$ is bounded and analytic for $\lambda\in S$ as a family of operators $\rho R_0(\lambda)\rho: L^2(\mathbb R^{2})\rightarrow L^2(\mathbb R^{2})$. It is easy to see that in the presence of the two cut-off functions $\rho$, the operator $\rho R_0(\lambda)\rho$ is bounded.
For the analyticity, it suffices to prove that the kernel $\rho H_0^{(1)}(\lambda |x-y|)\rho$ is analytic.
Since either $\Im\Big(\frac{t}{2} - {\rm i}\lambda |x - y|\Big)<0$ in \eqref{kernel}  for $\lambda\in \{\lambda: \arg\lambda\in [-\frac{\pi}{4}, \frac{\pi}{4}] \}$
or $\Im\Big(\frac{t}{2} - {\rm i}\lambda |x - y|\Big)>0$ for $\lambda\in \{\lambda: \arg\lambda\in [-\frac{3\pi}{4}, \frac{5\pi}{4}] \}$, 
one has that $\Big( \frac{t}{2} - {\rm i}\lambda |x - y| \Big)^{-\frac{1}{2}}$ is also analytic if one chooses the branch cut of $\sqrt{z}, z\in S$ as $\mathbb C\setminus [0, +\infty)$ fixed by the condition $\sqrt{1} = 1$. As a consequence, we have that
$\rho R_0(\lambda)\rho$ is analytic for $\lambda\in S$.

For $\Im\lambda>0$, as the term $e^{{\rm i} \lambda |  x -   y|}$ will be exponentially decaying for large $|  x-  y|$, it is easy to see from \eqref{kernel} that 
$\|R_0 (\lambda)\|_{L^2(\mathbb R^{2})\to L^2(\mathbb R^{2})} = \mathcal{O}(1/|\lambda|^{1/2})$.
Fix $\rho\in C_0^\infty(\mathbb R^{2})$. One has for $f\in L^2(\mathbb R^2)$ that
\begin{align}\label{kernel2}
(\rho R_0 (\lambda)\rho f)(  x) = \int_{\mathbb R^{2}}\rho(  x)R_0 (\lambda,   x,   y) \rho(  y)f(  y){\rm d}  y.
\end{align}
From the estimate \eqref{est} one has that $\rho R_0 (\lambda)\rho: L^2(\mathbb R^{2}) \to L^2(\mathbb R^{2})$ is bounded with the estimate
\begin{align*}
&\|\rho R_0(\lambda) \rho\|_{L^2(\mathbb R^{2})\rightarrow L^2(\mathbb R^{2})} \\
&\lesssim |\lambda|^{-\frac{1}{2}} e^{L(\Im\lambda)_-}\Big(\int_{\mathbb R^{2}}\int_{\mathbb R^{2}} \rho^2(  x)\frac{1}{|  x-  y|}\rho^2(  y){\rm d}  x{\rm d}  y\Big)^{1/2}\\
&\lesssim |\lambda|^{-\frac{1}{2}} e^{L(\Im\lambda)_-}.
\end{align*}
As the complex function $F(\lambda):=\langle \rho R_0(\lambda)\rho f, g\rangle_{L^2(\mathbb R^{2})}$ is analytic in $S$
for given $f, g\in L^2(\mathbb R^{2})$,
we have that $\rho R_0 (\lambda)\rho$ is an analytic family of bounded operators for $\lambda\in S$.

Next we prove \eqref{free}. We start with the case $j=2$. Choose $\tilde{\rho}\in C^\infty_0(\mathbb R^2)$ such that $\tilde{\rho} = 1$ near the support of $\rho$. 
By the standard elliptic estimate \cite[(7.13)]{stein} one has
\[
\|\rho u\|_{H^2(\mathbb R^2)}\leq C\big(\|\tilde{\rho}u\|_{L^2(\mathbb{R}^{2})} + \|\tilde{\rho}\Delta u\|_{L^2(\mathbb{R}^{2})}\big).
\]
Letting $u = R_0 (\lambda)(\rho f)$ in the above estimate one has 
\[
\|\rho R_0(\lambda)(\rho f)\|_{H^2(\mathbb R^2)}\leq C\big(\|\tilde{\rho}R_0 (\lambda)(\rho f)\|_{L^2(\mathbb R^{2})} + \|\tilde{\rho}\Delta(R_0 (\lambda)(\rho f))\|_{L^2(\mathbb R^{2})}\big).
\]
A direct calculation gives
\begin{align*}
\|\tilde{\rho}\Delta(R_0 (\lambda)(\rho f))\|_ {L^2(\mathbb R^2)}&=\| \rho f + \tilde{\rho}\lambda^2R_0 (\lambda)(\rho f)\|_{L^2(\mathbb R^2)}\\& \lesssim 
|\lambda|^{-\frac{1}{2}}(1+
|\lambda|^2) e^{L(\Im\lambda)_-}\|f\|_{L^2(\Omega)},
\end{align*}
which gives 
\[
\|\rho R_0 (\lambda)\rho\|_{L^2(\Omega)\rightarrow H^2(\Omega)}\lesssim 
|\lambda|^{-\frac{1}{2}}(1+ |\lambda|^2)
e^{L(\Im\lambda)_-}.
\]
The cases for $j = 1$ can be proved using the interpolation between $j=0$ and $j=2$. The proof is completed.

\end{proof}

Now we present the proof of Theorem \ref{meromorphic}, which provides a resonance-free region and resolvent estimates in this region for $R_V(\lambda)$.
The proof follows by utilizing the perturbation arguments in \cite{Dyatlov}.

We first present the analytic Fredholm theory below which will used in the subsequent analysis. The result is classical and the proof may be found in many references, e.g., \cite[Theorem 8.26]{CK}.

\begin{proposition}\label{AFT}
Let $D$ be a domain in $\mathbb C$ and let $K: D \rightarrow \mathcal{L}(X)$ be an operator
valued analytic function such that $K(z)$ is compact for each $z\in D$. Then either
\begin{itemize}

\item[(a)] $(I - K(z))^{-1}$ does not exist for any $z\in D$ or

\item[(b)] $(I - K(z))^{-1}$ exists for all $z\in D\backslash S$ where $S$ is a discrete subset of $D$.

\end{itemize}
Here $X$ is a Banach space and $\mathcal{L}(X)$ denotes the Banach space
of bounded linear operators mapping the Banach space $X$ into itself.
\end{proposition}

Let $\rho\in C_0^\infty(\mathbb R^2)$ with $\rho = 1$ on $\text{supp}V$. One has the following equality 
\begin{align}\label{equality}
(-\Delta + V(x) -\lambda^2)R_0(\lambda) = (-\Delta - \lambda^2) R_0(\lambda) + V(x)R_0(\lambda) = I + V(x)R_0(\lambda).
\end{align}
As for $\Im\lambda\gg 1$ 
\begin{align*}
\|VR_0(\lambda)\|_{L^2(\mathbb R^2)\rightarrow L^2(\mathbb R^2)} \leq \|V\|_{L^\infty(\mathbb R^2)}\| R_0(\lambda)\|_{L^2(\mathbb R^2)\rightarrow L^2(\mathbb R^2)}\leq\frac{\|V\|_{L^\infty(\mathbb R^2)}}
{\sqrt{|\lambda|}}\leq\frac{1}{2},
\end{align*}
using the Neumann series argument one has that the operator $I + VR_0(\lambda)$ is invertible with
\[
(I + VR_0(\lambda))^{-1} = \sum_{k=0}^\infty (-1)^k (VR_0(\lambda))^k.
\]
Thus, from \eqref{equality} one has that
\[
 R_V(\lambda) = R_0(\lambda) (I + VR_0(\lambda))^{-1}: L^2(\mathbb R^2)\rightarrow L^2(\mathbb R^2)
\]
are well-defined bounded operators for $\Im\lambda\gg 1$.

Define the following meromorphic family of operators for $\lambda\in S$:
\[
T(\lambda) = VR_0(\lambda): L_{\rm comp}^2(\mathbb R^2)\rightarrow L^2_{\rm comp}(\mathbb R^2).
\]
As $\rho T(\lambda) = \rho V R_0(\lambda) = VR_0(\lambda) = T(\lambda)$, one has $(1 - \rho)T(\lambda) = 0$, 
\[
(I + T(\lambda)(1 - \rho))^{-1} = I - T(\lambda)(1 - \rho)
\]
and
\[
(I + T(\lambda))^{-1} = (I + T(\lambda)\rho)^{-1} (I - T(\lambda)(1 - \rho)).
\]
Therefore for $\Im\lambda\gg 1$
\begin{align}\label{expression1}
R_V(\lambda) = R_0(\lambda)(I + T(\lambda))^{-1} =  R_0(\lambda) (I + T(\lambda)\rho)^{-1} (I - T(\lambda)(1 - \rho)).
\end{align}
Note that from Theorem \ref{free_estimate_2d}
\[
I - T(\lambda)(1 - \rho): L^2_{\rm comp}(\mathbb R^2)\rightarrow L^2_{\rm comp}(\mathbb R^2)
\]
and
\[
R_0(\lambda):  L^2_{\rm comp}(\mathbb R^2)\rightarrow H^2_{\rm loc}(\mathbb R^2)
\]
are both meromorphic operators for $\lambda \in S$. Hence in order to obtain the meromorphic continuation of $R_V(\lambda)$ as bounded operators
\[
R_V(\lambda): L^2_{\rm comp}(\mathbb R^2)\rightarrow L^2_{\rm loc}(\mathbb R^2)
\]
from $\{\lambda\in S: \Im\lambda>0\}$ to $S$, from the expression \eqref{expression1} of $R_V(\lambda)$ it suffices to prove
\[
(I + T(\lambda)\rho)^{-1} : L^2_{\rm comp}(\mathbb R^2)\rightarrow L^2_{\rm comp}(\mathbb R^2)
\]
is a meromorphic family of operators on $S$. Since by $V(  x) = V(   x)\rho(   x)$ we have $T(\lambda)\rho = V\rho R_0(\lambda)\rho$ and 
\begin{align*}
\|T(\lambda)\rho\|_{L^2_{\rm comp}(\mathbb R^2)\rightarrow L^2_{\rm comp}(\mathbb R^2)} &= \|V\rho R_0(\lambda)\rho\|_{L^2_{\rm comp}(\mathbb R^2)\rightarrow L^2_{\rm comp}(\mathbb R^2)}\\
&\leq\|V\|_{L^\infty} \|\rho R_0(\lambda)\rho\|_{L^2_{\rm comp}(\mathbb R^2)\rightarrow L^2_{\rm comp}(\mathbb R^2)}\\
&\leq C|\lambda|^{-1/2}\leq\frac{1}{2}
\end{align*}
for  $|\lambda|\gg 1$. Hence it follows from the Neumann series argument that the operator $(I + T(\lambda)\rho)^{-1}: L^2(\mathbb R^2)\rightarrow L^2(\mathbb R^2)$ exists for $|\lambda|\gg 1$. Moreover, for any $\lambda\in S$ the operator $T(\lambda)\rho = V\rho R_0(\lambda)\rho$ is compact on $L^2(\mathbb R^2)$ by the resolvent estimate \eqref{free}. Therefore, it follows from the analytic Fredholm theorem (cf Proposition \ref{AFT}) that $(I + T(\lambda)\rho)^{-1}: L^2(\mathbb R^2)\rightarrow L^2(\mathbb R^2)$ is meromorphic on $S$.

Finally, it remains to show that $(I + T(\lambda)\rho)^{-1}$ is 
 $L^2_{\rm comp}(\mathbb R^2)\rightarrow L^2_{\rm comp}(\mathbb R^2)$. In fact, we can choose $\chi, \tilde{\chi}\in C_{0}^\infty(\mathbb R^2)$ such that $\chi\rho = \rho$ and $\tilde{\chi}\chi = \chi$, then $(1 - \tilde{\chi})\rho = 0$. Moreover, when $|\lambda|\gg 1$, by the Neumann series argument and $V\rho = V$, we have
\begin{align}\label{cutoff}
(1 - \tilde{\chi})(I + T(\lambda)\rho)^{-1}\chi &= (1 - \tilde{\chi})\chi + \sum_{k=1}^\infty (-1)^k (1 - \tilde{\chi}) (T(\lambda)\rho)^k\chi\notag\\
&= \sum_{k=1}^\infty (-1)^k (1 - \tilde{\chi}) (V\rho R_0(\lambda)\rho)^k\chi\notag\\
&= \sum_{k=1}^\infty (-1)^k (1 - \tilde{\chi}) (V\rho R_0(\lambda)\rho)(V\rho R_0(\lambda)\rho)^{k-1}\chi\notag\\
&=0,
\end{align}
where the last equality uses $(1 - \tilde{\chi})\rho = 0$. By the analytic continuation, \eqref{cutoff} remains true for all $\lambda$ where $T(\lambda)$ is analytic. Therefore, by the expression \eqref{expression1} of $R_V(\lambda)$ we obtain that $R_V(\lambda)$ is meromorphic for $\lambda\in S$ as a family of operators $ L^2_{\rm comp}(\mathbb R^2)\rightarrow H^2_{\rm loc}(\mathbb R^2)$.

Now we prove the resolvent estimates. For $\lambda\in \Omega_\delta\cap S$, one has
\begin{align*}
\|V R_0(\lambda)\rho\|_{L^2(\mathbb R^2)\rightarrow L^2(\mathbb R^2)} &= \|V\rho R_0(\lambda)\rho_1\|_{L^2(\mathbb R^2)\rightarrow L^2(\mathbb R^2)}\\
&\lesssim \|V\|_{L^\infty}|\lambda|^{-1/2} e^{L(\Im\lambda)_-}\\
&\lesssim \|V\|_{L^\infty}|\lambda|^{-1/2}e^{L\delta\text{log}|\lambda|}\\
&\lesssim  \|V\|_{L^\infty}|\lambda|^{-1/4}\leq\frac{1}{2},
\end{align*}
where we let $C_0\gg 1$ and $\delta<\frac{1}{4L}.$ 
Hence by the Neumann series argument we can prove that the inverse operator $(I + VR_0(\lambda)\rho)^{-1}$ exists for all $\lambda\in\Omega_\delta\cap S$, and
\begin{align}\label{free4}
\|(I + VR_0(\lambda)\rho)^{-1}\|_{L_{\rm comp}^2(\mathbb R^2)\rightarrow L_{\rm comp}^2(\mathbb R^2)} = \|(I + V\rho R_0(\lambda)\rho)^{-1}\|_{L_{\rm comp}^2(\mathbb R^2)\rightarrow L_{\rm comp}^2(\mathbb R^2)}\leq 2.
\end{align}
One has
\[
\rho R_V(\lambda) \rho = \rho R_0(\lambda) \rho (I + VR_0(\lambda)\rho_1)^{-1} (I - VR_0(\lambda)(1 -\rho_1))\rho.
\]
Furthermore, by the estimate \eqref{free} we obtain 
\begin{align}\label{free3}
\|\rho R_0(\lambda)\rho\|_{L^2(\mathbb R^2)\to L^2(\mathbb R^2)} \leq C |\lambda|^{-1/2} e^{L(\Im\lambda)_-}.
\end{align}
Combining  \eqref{free3} and \eqref{free4} one has the desired estimate for $j = 0$ as
\[
\|\rho R_V(\lambda)\rho\|_{L^2(\mathbb R^2)\rightarrow L^2(\mathbb R^2)} \leq C |\lambda|^{-1/2} e^{L(\Im\lambda)_-}.
\]
 For the case $j = 2$, let $\tilde{\rho}\in C_{0}^\infty(\mathbb R^2)$ such that $\tilde{\rho} = 1$ on $\text{supp}\rho$. One has 
 from the standard elliptic estimate \cite[(7.13)]{stein} that
\begin{align*}
\|\rho R_V(\rho f)\|_{H^2(\mathbb R^2)} &\lesssim  \|\tilde{\rho} R_V(\lambda)(\rho f)\|_{L^2(\mathbb R^2)} + \|\tilde{\rho} \Delta R_V(\lambda)
(\rho f)\|_{L^2(\mathbb R^2)}\\
&\lesssim \|\tilde{\rho} R_V(\lambda)(\rho f)\|_{L^2(\mathbb R^2)} + \|\tilde{\rho} (-\Delta + V) R_V(\lambda)(\rho f)\|_{L^2(\mathbb R^2)} + \|\tilde{\rho} V R_V(\lambda)(\rho f)\|_{L^2(\mathbb R^2)}\\
&\lesssim \|\tilde{\rho} R_V(\lambda)(\rho f)\|_{L^2(\mathbb R^2)} + |\lambda|^2 \|\tilde{\rho} R_V(\lambda)(\rho f)\|_{L^2(\Omega)} 
+ \|\rho f\|_{L^2(\mathbb R^2)}\\ 
&\quad + \|\tilde{\rho} V R_V(\lambda)(\rho f)\|_{L^2(\mathbb R^2)}\\
&\lesssim  (1+ |\lambda|)^2  \|\tilde{\rho} R_V(\lambda)(\rho f)\|_{L^2(\mathbb R^2)} \\
&\lesssim (1+ |\lambda|)^2 |\lambda|^{-1/2} e^{L(\Im\lambda)_-}\|f\|_{L^2(\mathbb R^2)}
\end{align*}
for $\lambda\in\Omega_\delta\cap S$. Finally, the cases of $j =1$ follows by an application of the interpolation between $j = 0$ and $j = 2$. 
The proof is completed.

\end{document}